\documentclass[draft]{amsart}
\usepackage{amsmath,amssymb,amscd}
\newtheorem{formula}{}[section]
\newtheorem{proposition}[formula]{Proposition}
\newtheorem{corollary}[formula]{Corollary}
\newtheorem{lemma}[formula]{Lemma}
\newtheorem{theorem}[formula]{Theorem}
\theoremstyle{definition}
\newtheorem{definition}[formula]{Definition}
\newtheorem{example}[formula]{Example}
\theoremstyle{remark}
\newtheorem*{remark}{Remark}

\newcommand{\Tor}{\mathop{\rm Tor}\nolimits}
\newcommand{\Krull}{\mathop{\rm Krull}}
\newcommand{\bideg}{\mathop{\rm bideg}}

\renewcommand{\l}{\lambda}
\renewcommand{\emptyset}{\varnothing}
\newcommand{\C}{\mathbb C}
\newcommand{\R}{\mathbb R}
\newcommand{\Z}{\mathbb Z}
\newcommand{\D}{\Delta}
\newcommand{\G}{\Gamma}
\newcommand{\ZP}{\mathcal Z_P}
\newcommand{\T}{\mathcal T}

\renewcommand{\>}{\rangle}

\newcommand{\ep}{\nolinebreak\quad$\square$}
\newcommand{\eqn}{\refstepcounter{equation}\leqno{(\arabic{equation})}}

\begin{document}

\title{Torus actions and combinatorics of polytopes}
\author{Victor M. Buchstaber}
\author{Taras E. Panov}
\thanks{Partially supported by
the Russian Foundation for Fundamental Research, grant no. 99-01-00090.}
\address{Department of Mathematics and Mechanics, Moscow
State University, 119899 Moscow, RUSSIA}
\email{tpanov@mech.math.msu.su \quad buchstab@mech.math.msu.su}

\begin{abstract}
An $n$-dimensional polytope $P^n$ is called {\it simple} if exactly
$n$ codimension-one faces meet at each vertex.
The lattice of faces of a simple polytope $P^n$ with $m$
codimension-one faces defines an arrangement of
even-dimensional planes in
$\R^{2m}$. We construct a {\it free} action of the group $\R^{m-n}$
on the complement of this arrangement. The corresponding
quotient is a smooth manifold $\ZP$ invested with a canonical action of the
compact torus $T^m$ with the orbit space $P^n$.  For each smooth projective
{\it toric variety} $M^{2n}$ defined by a simple polytope $P^n$ with the given
lattice of faces there exists a subgroup $T^{m-n}\subset T^m$ acting freely
on $\ZP$ such that $\ZP/T^{m-n}=M^{2n}$.  We calculate the cohomology ring of
$\ZP$ and show that it is isomorphic to the cohomology ring of the {\it face
ring} of $P^n$ regarded as a module over the polynomial ring. In this way the
cohomology of $\ZP$ acquires a {\it bigraded} algebra structure, and the
additional grading allows to catch the combinatorial invariants of the
polytope. At the same time this gives an example of explicit calculation of
the cohomology ring for the complement of an arrangement of planes,
which is of independent interest.
\end{abstract}

\maketitle

\section*{Introduction}

In this paper we study relations between the
algebraic topology of manifolds and the combinatorics of polytopes.
Originally, this research was inspired by the results of the toric variety
theory. The main object of our study is the smooth manifold defined by the
combinatorial structure of a simple polytope. This manifold is equipped with
a natural action of the compact torus $T^m$.

We define an $n$-dimensional {\it convex polytope} as a bounded set in $\R^n$
that is obtained as the intersection of a finite number of half-spaces. So,
any convex polytope is bounded by a finite number of hyperplanes. A convex
$n$-dimensional polytope is called {\it simple} if there exactly $n$
co\-di\-men\-si\-on-one faces (or {\it facets}) meet at each vertex. The
bounding hyperplanes of a simple polytope are in general position at each
vertex. A convex polytope can be also defined as the convex hull of a set of
points in $\R^n$. If these points are in general position, the resulting
polytope is called {\it simplicial}, since all its faces are simplices. For
each simple polytope there is defined the {\it dual} (or {\it polar})
simplicial polytope (see Definition~\ref{dp}). The boundary of a simplicial
polytope defines a simplicial subdivision (triangulation) of a sphere.

We associate to each simple polytope $P^n$ with $m$ facets a smooth
$(m+n)$-dimensional manifold $\ZP$ with a canonical action of the compact
torus $T^m$. A~number of manifolds that play the important role in different
aspects of topology, algebraic and symplectic geometry appear as special
cases of the above manifolds $\ZP$ or as the quotients $\ZP/T^k$ for torus
subgroups $T^k\subset T^m$ acting on $\ZP$ freely. It turns out that no torus
subgroup of the rank $>m-n$ can act on $\ZP$ freely. We call the quotients of
$\ZP$ by tori of the maximal possible rank $m-n$ {\it quasitoric manifolds}.
The name refers to the fact that the important class of algebraic varieties
known to algebraic geometers as {\it toric manifolds} fits the above picture.
More precisely, one can use the above construction (i.e. the quotient of
$\ZP$ by a torus subgroup) to produce all smooth projective {\it toric
varieties} (see~[Da]), which we refer to as {\it toric manifolds}. The action
of $T^m$ on $\ZP$ induces an action of $T^n$ on the (quasi)toric manifold
$M^{2n}:=\ZP/T^{m-n}$ whit the same orbit space $P^n$. However, there are
combinatorial types of simple polytopes that can not be realized as the orbit
space for a quasitoric (and also toric) manifold. This means exactly that for
such combinatorial type $P$ it is impossible to find a torus subgroup
$T^{m-n}\subset T^m$ of rank $m-n$ that acts on the corresponding manifold
$\ZP$ freely. If the manifold $\ZP$ defined by a combinatorial simple
polytope $P$ admits a free action of a torus subgroup of rank $m-n$, then
different such subgroups may produce different quasitoric manifolds over
$P^n$, and some of them may turn out to be toric manifolds. Our quasitoric
manifolds originally appeared under the name ``toric manifolds" in~[DJ],
where the authors developed combinatorial and topological methods for
studying the corresponding torus actions. We use some results of~\cite{DJ} in
our paper.

Our approach to constructing manifolds $\ZP$ defined by simple polytopes is
based on one construction from algebraic geometry, which was used in~[Ba] for
studying toric varieties. Namely, the lattice of faces of a simple polytope
$P^n$ defines a certain affine algebraic set $U(P^n)\subset\C^m$ with an
action of the algebraic torus $(\C^{*})^m$. This set $U(P^n)$ is the
complement of a certain arrangement of planes in $\C^m$ defined by the
combinatorics of $P^n$.  Toric manifolds arise when one can find a subgroup
$C\subset(\C^{*})^m$ isomorphic to $(\C^{*})^{m-n}$ that acts on $U(P^n)$
freely. The crucial fact in our approach is that it is {\it always} possible
to find a subgroup $R\subset(\C^{*})^m$ isomorphic to $\R^{m-n}$ and acting
freely on $U(P^n)$. In this case the quotient manifold is defined, which we
refer to as the manifold defined by the simple polytope $P^n$. There is a
canonical action of the torus $T^m$ on this manifold, namely, that induced by
the standard action of $T^m$ on $\C^m$ by diagonal matrices.  The another
approach to construct $\ZP$ was proposed in~[DJ], where this manifold was
defined as the quotient space $\ZP=T^m\times P^n/\sim$ for a certain
equivalence relation $\sim$. We construct an equivariant embedding $i_e$ of
this manifold into $U(P^n)\subset\C^m$ and show that for the above subgroup
$R\cong\R^{m-n}$ the composite map $\ZP\to U(P^n)\to U(P^n)/R$ of the
embedding and the orbit map is a homeomorphism. Hence, from the topological
viewpoint, both approaches produce the same manifold. This is what we refer
to as the manifold defined by simple polytope $P^n$ and denote by $\ZP$
throughout our paper.

The analysis of the above constructions shows that we can replace the
$m$-dimensional complex space $\C^m\cong(\R^2)^m$ by a space $(\R^k)^m$
with arbitrary $k$.  Indeed, we may construct the open subset
$U_{(k)}(P^n)\subset(\R^k)^m$ determined by the lattice of faces of $P^n$ in
the same way as in the case of $\C^m$ (i.e. $U_{(k)}(P^n)$ is the complement
of a certain set of planes, see Definition~\ref{defU}). The multiplicative
group $(\R_>)^m\cong\R^m$ acts on $(\R^k)^m$ diagonally (i.e. as the product
of $m$ standard actions of $\R_>$ on $\R^k$ by dilations). For this action it
is also possible to find a subgroup $R\subset(\R_>)^m$ isomorphic to
$\R^{m-n}$ that acts on $U_{(k)}(P^n)$ freely. The corresponding quotient
$U_{(k)}(P^n)/R$ is now of dimension $(k-1)m+n$ and is invested with an
action of the group $O(k)^m$ (the product of $m$ copies of the orthogonal
group). This action is induced by the diagonal action of $O(k)^m$ on
$(\R^k)^m$. In the case $k=2$ the above considered action of the torus $T^m$
on $U(P^n)/R$ is just the action of $SO(2)^m\subset O(2)^m$. In the case
$k=1$ we obtain for any simple polytope $P^n$ a smooth $n$-dimensional
manifold $\mathcal Z^n$ with an action of the group $(\Z/2)^m$ whose orbit
space is $P^n$. This manifold is known as the universal Abelian cover of
$P^n$ regarded as a right-angled Coxeter orbifold (or manifold with corners).
The analogues of quasitoric manifolds in the case $k=1$ are the so-called
{\it small covers}. A small cover is a $n$-dimensional manifold $M^n$ with
an action of $(\Z/2)^n$ whose orbit space is $P^n$. The name refers to the
fact that any cover of $P^n$ by a manifold must have at least $2^n$ sheets.
The case $k=1$ was detailedly treated in~\cite{DJ}, along with quasitoric
manifolds. The another case of particular interest is $k=4$, since $\R^4$
can be viewed as a one-dimensional quaternionic space. In this paper we study
the case $k=2$, which all the constructions below refer to.

One of our main goals here is to study relations between the combinatorics of
simple polytopes and the topology of the above described manifolds. There is
a well-known important algebraic invariant of a simple polytope: a graded
ring $k(P)$ (here $k$ is any field), called {\it the face ring} (or the {\it
Stanley--Reisner ring}), see~\cite{St}. This is the quotient of the
polynomial ring $k[v_1,\ldots,v_m]$ by a homogeneous ideal determined by the
lattice of faces of a polytope (see Definition~\ref{frpol}). The cohomology
modules $\Tor^{-i}_{k[v_1,\ldots,v_m]}\bigl(k(P),k\bigr)$, $i>0$, are of
great interest to algebraic combinatorists.  Some results on the
corresponding Betti numbers $\beta^i\bigl(k(P)\bigr)=
\dim_k\Tor^{-i}_{k[v_1,\ldots,v_m]}\bigl(k(P),k\bigr)$ can be found in~[St].
We show that the bigraded $k$-module
$\Tor_{k[v_1,\ldots,v_m]}\bigl(k(P),k\bigr)$ can be endowed with a bigraded
$k$-algebra structure and its totalized graded algebra is isomorphic to the
cohomology algebra of $\ZP$. Therefore, the cohomology of $\ZP$ possesses a
canonical {\it bigraded} algebra structure. The proof of these facts uses the
Eilenberg--Moore spectral sequence. This spectral sequence is usually applied
in algebraic topology as a powerful tool for calculating the cohomology of
homogeneous spaces for Lie group actions (see e.g.,~[Sm]). So, it was
interesting for us to discover a quite different application of this spectral
sequence. In our situation the $E_2$ term of the spectral sequence is exactly
$\Tor_{k[v_1,\ldots,v_m]}\bigl(k(P),k\bigr)$, and the spectral sequence
collapses in the $E_2$ term. Using the Koszul complex as a resolution while
calculating the $E_2$ term, we show that the above bigraded algebra is the
cohomology algebra of a certain bigraded complex defined in purely
combinatorial terms of the polytope $P^n$ (see Theorem~\ref{mult}).
Therefore, our bigraded cohomology algebra of $\ZP$ contains all the
combinatorial data of $P^n$. In particular, it turns out that the well-known
Dehn--Sommerville equations for a simple polytope $P^n$ follow directly from
the bigraded Poincar\'e duality for $\ZP$. Given the corresponding bigraded
Betti numbers one can compute the numbers of faces of $P$ of fixed dimension
(the so-called {\it $f$-vector} of the polytope). Many combinatorial results,
such as the Upper Bound for the number of faces of a simple polytope, can be
interpreted nicely in terms of the cohomology of the manifold $\ZP$.

Moreover, since the homotopy equivalence $\ZP\simeq U(P^n)$ holds, our
calculation of the cohomology is also applicable to the set $U(P^n)$. As it
was mentioned above, $U(P^n)$ is the complement of an arrangement of planes
in $\C^m$ defined by the combinatorics of $P^n$. Hence, here we have a
special case of the well-known general problem of calculating the cohomology
of the complement of an arrangement of planes. In~[GM,~part~III] the
corresponding Betti numbers were calculated in terms of the cohomology of a
certain simplicial complex.  In our case special properties of the
arrangement defined by a simple polytope allow to calculate the {\it
cohomology ring} of the corresponding complement much more explicitly.

Problems considered here were discussed in the first author's talk on the
conference ``Solitons, Geometry and Topology" devoted to the jubilee of our
Teacher Sergey Novikov. A part of the results of this paper were announced
in~\cite{BP}.

\section{Main constructions and definitions}

\subsection{Simple polytopes and their face rings.}

Let $P^n$ be a simple polytope and let $f_i$ be its number of
codimension $(i+1)$ faces, $0\le i\le n-1$. We refer to the integer vector
$(f_0,\ldots,f_{n-1})$ as the {\it $f$-vector} of $P^n$. It is convenient
to put $f_{-1}=1$. Along with the $f$-vector we also consider the {\it
$h$-vector} $(h_0,\ldots,h_n)$ whose components $h_i$ are
defined from the equation
\begin{equation}
\label{hvector}
  h_0t^n+\ldots+h_{n-1}t+h_n=(t-1)^n+f_0(t-1)^{n-1}+\ldots+f_{n-1}.
\end{equation}
Therefore, we have
\begin{equation}
\label{hf}
  h_k=\sum_{i=0}^k(-1)^{k-i}\binom{n-i}{k-i}f_{i-1}.
\end{equation}

We fix a commutative ring $k$, which we refer to as the ground
ring. A graded ring called {\it face ring} is associated to
$P^n$. More precisely, let
${\mathcal F}=(F_1,\ldots,F_m)$ be the set of codimension-one faces of $P^n$,
$m=f_0$. Form the polynomial ring $k[v_1,\ldots,v_m]$, where the $v_i$ are
regarded as indeterminates corresponding to the facets $F_i$.

\begin{definition}
\label{frpol}
  The {\it face ring\/} $k(P)$ of a simple polytope $P$ is defined to be the
  ring $k[v_1,\ldots,v_m]/I$, where
  $$
    I=\left(v_{i_1}\ldots v_{i_s}:\;i_1<i_2<\ldots <i_s,\;
    F_{i_1}\cap F_{i_2}\cap\cdots\cap F_{i_s}=\emptyset\right).
  $$
\end{definition}
Note that the face ring is determined only by a combinatorial type of $P^n$
(i.e. by its lattice of faces).

In the literature (see \cite{St}), the face ring is usually defined for
simplicial complexes, as follows. Let $K$ be a finite simplicial complex with
the vertex set $\{v_1,\ldots,v_m\}$. Form a polynomial ring
$k[v_1,\ldots,v_m]$ where the $v_i$ are regarded as indeterminates.

\begin{definition}
\label{frsim}
  The {\it face ring\/} $k(K)$ of a simplicial complex $K$ is
  the quotient ring $k[v_1,\ldots,v_m]/I$, where
  $$
    I=\left(v_{i_1}\ldots v_{i_s}:\;i_1<i_2<\ldots <i_s,\;
    \{v_{i_1},\ldots,v_{i_s}\}\; \text{ does not span a simplex in }K
    \right).
  $$
\end{definition}
We regard the indeterminates $v_i$ in
$k[v_1,\ldots,v_m]$ as being of degree two; in this way $k(K)$, as well
as $k(P)$, becomes a graded ring.

\begin{definition}
\label{dp}
  Given a convex polytope $P^n\subset\R^n$, the
  {\it dual} (or {\it polar}) polytope
  $(P^n)^{*}\subset(\R^n)^{*}$ is defined as follows
  $$
    (P^n)^{*}=\{x'\in(\R^n)^{*}\::\:\;\langle x',x\rangle\le1
    \:\text{ for all }\:x\in P^n\}.
  $$
\end{definition}
It can be shown (see~\cite{Br}) that the above set is indeed a convex
polytope. In the case of simple $P^n$ the dual polytope $(P^n)^{*}$ is
simplicial and its $i$-dimensional faces (simplices)
are in one-to-one correspondence with the faces of $P^n$ of codimension $i+1$.
The boundary complex of $(P^n)^{*}$ defines a simplicial subdivision
(triangulation) of a $(n-1)$-dimensional sphere $S^{n-1}$, which we denote
$K_P$. In this situation both Definitions~\ref{frpol} and~\ref{frsim}
provide the same ring: $k(P)=k(K_P)$. The face rings of
simple polytopes have very specific algebraic properties. To
describe them we need some commutative algebra.

Now suppose that $k$ is a field and let $R$ be a graded algebra over $k$.
Let $n$ be the maximal number of algebraically independent
elements of $R$ (this number is known as the {\it Krull dimension} of $R$,
denoted $\Krull R$). A sequence $(\l_1,\ldots,\l_k)$ of homogeneous elements
of $R$ is called a {\it regular sequence}, if $\l_{i+1}$ is not a zero
divisor in $R/(\l_1,\ldots,\l_i)$ for each $i$ (in the other words, the
multiplication by $\l_{i+1}$ is a monomorphism of $R/(\l_1,\ldots,\l_i)$
into itself). It can be proved that $(\l_1,\ldots,\l_k)$ is a regular
sequence if and only if $\l_1,\ldots,\l_k$ are algebraically independent and
$R$ is a free $k[\l_1,\ldots,\l_n]$-module. The
notion of regular sequence is of great importance for algebraic
topologists (see, for instance,~\cite{La},~\cite{Sm}).
A sequence $(\l_1,\ldots,\l_n)$ of homogeneous elements
of $R$ is called a {\it homogeneous system of parameters} (hsop), if the
Krull dimension of $R/(\l_1,\ldots,\l_n)$ is zero. The $k$-algebra $R$ is
{\it Cohen--Macaulay} if it admits a regular sequence $(\l_1,\ldots,\l_n)$ of
$n=\Krull R$ elements (which is then automatically a hsop). It follows from
the above that $R$ is Cohen--Macaulay if and only if there exists a sequence
$(\l_1,\ldots,\l_n)$ of algebraically independent homogeneous elements of $R$
such that $R$ is a finite-dimensional free $k[\l_1,\ldots,\l_n]$-module.

In our case the following statement holds (see \cite{St}).
\begin{proposition}
  The face ring $k(P^n)$ of a simple polytope $P^n$ is a Cohen--Macaulay
  ring.\ep
\end{proposition}

In what follows we need two successive generalizations of a simple polytope.
As it was mentioned in the introduction, the bounding hyperplanes of a simple
polytope are in general position at each vertex.  First, we define a {\it
simple polyhedron} as a convex set in $\R^n$ (not necessarily bounded)
obtained as the intersection of a finite number of half-spaces with the
additional condition that no more than $(n+1)$ hyperplanes intersect in one
point. The faces of a simple polyhedron are defined obviously; all of them
are simple polyhedra as well. It is also possible to define the
$(n-1)$-dimensional simplicial complex $K_P$ dual to the boundary of a simple
polyhedron $P^n$. (And again the $i$-dimensional simplices of $K_P$ are in
one-to-one correspondence with the faces of $P^n$ of codimension $i+1$.)
However, the simplicial complex $K_P$ obtained in such way not necessarily
defines a triangulation of a $(n-1)$-dimensional sphere $S^{n-1}$.

\begin{example}
  The simple polyhedron
  $$
    {\mathbb R}^n_+=\{(x_1,\ldots,x_n)\in{\mathbb R}^n:\;x_i\ge0\},
  $$
  bounded by $n$ coordinate hyperplanes will appear many times
  throughout our paper. The simplicial complex dual to its boundary is
  a $(n-1)$-dimensional simplex $\D^{n-1}$.
\end{example}

We note, however, that not any $(n-1)$-dimensional simplicial complex
can be obtained as the dual to a $n$-dimensional simple polyhedron. Because
of this, we still need to generalize the notion of a simple polytope (and
simple polyhedron). In this way we come to the notion of a {\it simple
polyhedral complex}. Informally, a simple polyhedral complex of dimension $n$
is ``the dual to a general $(n-1)$-dimensional simplicial complex". We take
its construction from~\cite{DJ}.  Let $K$ be a simplicial complex of
dimension $n-1$ and let $K'$ be its barycentric subdivision. Hence, the
vertices of $K'$ are simplices $\D$ of the complex $K$, and the simplices of
$K'$ are sets $(\Delta_1,\Delta_2,\ldots,\Delta_k)$, $\Delta_i\in K$, such
that $\Delta_1\subset\Delta_2\subset\ldots\subset\Delta_k$. For each simplex
$\Delta\in K$ denote by $F_{\Delta}$ the subcomplex of
$K'$ consisting of all simplices of $K'$ of the form
$\Delta=\Delta_0\subset\Delta_1\subset\ldots\subset\Delta_k$. If
$\Delta$ is a $(k-1)$-dimensional simplex, then we refer to
$F_{\Delta}$ as a face of codimension $k$. Let $P_K$ be the cone over
$K$. Then this $P_K$ together with its decomposition into
``faces" $\{F_{\Delta}\}_{\Delta\in K}$ is said to be a {\it simple
polyhedral complex}. Any simple polytope $P^n$ (as well as a simple
polyhedron) can be obtained by applying this construction to the
simplicial complex $K^{n-1}$ dual to the boundary $\partial P^n$.

\subsection{The topological spaces defined by simple polytopes.}

Following \cite{DJ}, in this subsection we associate two topological spaces,
$\ZP$ and $B_TP$, to any simple polyhedral complex $P$ (in particular, to any
simple polytope).

Let $T^m=S^1\times\ldots\times S^1$ be the $m$-dimensional compact
torus. Let ${\mathcal F}=(F_1,\ldots,F_m)$ denote, as before, the set of
codimension-one faces of $P^n$ (or the vertex set of the
dual simplicial complex $K^{n-1}$). We consider the lattice $\Z^{m}$ of
one-parameter subgroups of $T^m$ and fix a one-to-one correspondence between
the facets of $P^n$ and the elements of a basis
$\{e_1,\ldots,e_m\}$ in $\Z^{m}$. Now we can define the canonical coordinate
subgroups $T^k_{i_1,\ldots i_k}\subset T^m$ as the tori corresponding to the
coordinate subgroups of $\Z^{m}$ (i.e. to the subgroups spanned by
basis vectors $e_{i_1},\ldots,e_{i_k}$).

\begin{definition}
\label{defzp}
The equivalence relation $\sim$ on $T^m\times P^n$ is defined as follows
\begin{align*}
  &(g,p)\sim (h,q) \text{ iff }p=q\text{ and }g^{-1}h\in
  T^k_{i_1,\ldots,i_k},\\
  &\text{where $p$ lies in the relative interior of the face
  $F_{i_1}\cap\cdots\cap F_{i_k}$.}
\end{align*}
We associate to any simple polytope $P^n$ a topological space
$\ZP=(T^m\times P^n)/\!\sim$
\end{definition}
As it follows from the
definition, $\dim {\ZP}=m+n$ and the action of the torus $T^m$ on $T^m\times
P^n$ descends to an action of $T^m$ on $\ZP$.  In the case of simple
polytopes, the orbit space for this action is a $n$-dimensional ball invested
with the combinatorial structure of the polytope $P^n$ as described by the
following proposition.

\begin{proposition}
\label{combstr}
  Suppose that $P^n$ is a simple polytope. Then the action of $T^m$ on
  $\ZP$ has the following properties:
  \begin{enumerate}
  \item The isotropy subgroup of any point of $\ZP$ is a
    coordinate subgroup of $T^m$ of dimension $\le n$.
  \item The isotropy subgroups define the combinatorial structure of the
    polytope $P^n$ on the orbit space. More precisely, the interior of a
    codimension-$k$ face consists of orbits with the same $k$-dimensional
    isotropy subgroup. In particular, the action is free
    over the interior of the polytope.
  \end{enumerate}
\end{proposition}
\begin{proof}
This follows easily from the definition of $\ZP$.
\end{proof}
\medskip

Now we return to the general case of a simple polyhedral complex $P^n$. Let
$ET^m$ be the contractible space of the universal principal $T^m$-bundle over
$BT^m=(\C P^{\infty})^m$. Applying the Borel construction to the
$T^m$-space $\ZP$, we come to the following definition.

\begin{definition}
  The space $B_TP$ is defined as
  \begin{equation}
  \label{bun}
    B_TP=ET^m\times_{T^m}{\ZP}.
  \end{equation}
\end{definition}
Hence, the $B_TP$ is the total space of the bundle (with the fibre
$\ZP$) associated to the universal bundle via the action of
$T^m$ on $\ZP$. As it follows from the definition, the homotopy type of
$B_TP$ is determined by a simple polyhedral complex $P^n$.

\subsection{Toric and quasitoric manifolds.}

In the previous subsection we defined for any simple polytope $P^n$ a
space $\ZP$ with an action of $T^m$ and the combinatorial structure of $P^n$
in the orbit space (Proposition~\ref{combstr}). As we shall see in the next
section, this $\ZP$ turns out to be a smooth manifold. Another class of
manifolds possessing the above properties is well known in algebraic
geometry as toric manifolds (or non-singular projective toric varieties).
Below we give a brief review of them. The detailed background material on
this subject can be found in~\cite{Da},~\cite{Fu}.

\begin{definition}
  A {\it toric variety} is a normal algebraic variety $M$ containing the
  $n$-dimensional algebraic torus $(\C^{*})^n$ as a Zariski open subvariety,
  with the additional condition that the diagonal action of $(\C^{*})^n$ on
  itself extends to an action on the whole $M$ (so, the torus $(\C^{*})^n$ is
  contained in $M$ as a dense orbit).
\end{definition}

On any non-singular projective toric variety there exists a very ample line
bundle whose zero cohomology (the space of global sections) is
generated by the sections corresponding to the points with integer
coordinates inside a certain simple polytope with vertices in the integer
lattice $\Z^n\subset\R^n$. Conversely, there is an algebraic construction
which allows to produce a projective toric variety $M^{2n}$ of real
dimension $2n$ starting from a simple polytope $P^n$ with
vertices in $\Z^n$ (see, e.g.,~\cite{Fu}). However, the resulting variety
$M^{2n}$ is not necessarily non-singular. Namely, this construction gives
us a non-singular variety if
for each vertex of $P^n$ the normal covectors of $n$ facets meeting at this
vertex form a basis of the dual lattice $(\Z^n)^{*}$. A~toric variety is
not uniquely determined by the combinatorial type of a polytope: it depends
also on integer coordinates of vertices. Thus, a given {\it
combinatorial type} of a simple polytope (i.e. a lattice of faces) gives rise
to a number of toric varieties, one for each {\it geometrical} realization
with integer vertices. Some of these toric varieties may be non-singular.
However, there are combinatorial simple polytopes that produce only singular
toric varieties. The corresponding examples will be discussed below.

The algebraic torus contains the compact torus $T^n\subset(\C^{*})^n$, which
acts on a toric manifold as well. It can be proved that all isotropy
subgroups for this action are tori $T^k\subset T^m$ and the orbit space has
the combinatorial structure of $P^n$, as described
in the second part of Proposition~\ref{combstr} (here $P^n$ is the polytope
defined by the toric manifold as described above).
The action of $T^n$ on $M^{2n}$ is {\it locally equivalent to the
standard action of $T^n$ on $\C^n$} (by the diagonal matrices) in the
following sense: every point $x\in M^{2n}$ lies in some $T^n$-invariant
neighbourhood $U\subset M^{2n}$ which is $T^n$-equivariantly
homeomorphic to a certain ($T^n$-invariant) open subset $V\subset\C^n$.
Furthermore, there exists an explicit map $M^{2n}\to\R^n$ (the
{\it moment map}), with image $P^n$ and $T^n$-orbits as fibres
(see~\cite{Fu}). A toric manifold $M^{2n}$ (regarded as a smooth
manifold) can be obtained as the quotient space $T^n\times P^n/\!\sim$ for
some equivalence relation $\sim$ (see~\cite{DJ}; compare this with
Definition~\ref{defzp} of the space $\ZP$). Now, if we are interested
only in topological and combinatorial properties, then we should
not restrict ourselves to algebraic varieties; in this way, forgetting all
the algebraic geometry of $M^{2n}$ and the action of the algebraic torus
$(\C^{*})^n$, we come to the following definition.

\begin{definition}
  A {\it quasitoric} (or {topologically toric}) manifold over a simple
  polytope $P^n$ is a real orientable $2n$-dimensional manifold $M^{2n}$ with
  an action of the compact torus $T^n$ that is locally isomorphic to the
  standard action of $T^n$ on $\C^n$ and whose orbit space has the
  combinatorial structure of $P^n$ (in the sense of the second part of
  Proposition~\ref{combstr})
\end{definition}

Quasitoric manifolds were firstly introduced in~\cite{DJ} under the name
``toric manifolds". As it follows from the above discussion, all algebraic
non-singular toric varieties are quasitoric manifolds as well.  The converse
is not true: the corresponding examples can be found in~\cite{DJ}. One of the
most important result on quasitoric manifolds obtained there is the
description of their cohomology rings. This result generalizes the well-known
Danilov--Jurkiewicz theorem, which describes the cohomology ring of a
non-singular projective toric variety (see~\cite{Da}). In section~3 we give a
new proof of the result of~\cite{DJ} by means of the Eilenberg--Moore
spectral sequence (see Theorem~\ref{cohomM}). In the rest of this subsection
we describe briefly main constructions with quasitoric manifolds. The proofs
can be found in~\cite{DJ}.

\medskip

Suppose $M^{2n}$ is a quasitoric manifold over a simple polytope $P^n$, and
$\pi:M^{2n}\to P^n$ is the orbit map.
Let $F^{n-1}$ be a codimension-one face of $P^n$; then for any
$x\in\pi^{-1}({\mathop {\rm int}}F^{n-1})$ the isotropy subgroup at $x$ is
an independent of the choice of $x$ rank-one subgroup $G_F\in T^n$. This
subgroup is determined by a primitive vector $v\in\Z^n$. In this way we
construct a function $\l$ from the set $\mathcal F$ of codimension-one faces
of $P^n$ to primitive vectors in $\Z^n$.

\begin{definition}
\label{chf}
  The defined above function $\l:{\mathcal F}\to\Z^n$ is called the
  {\it characteristic function} of $M^{2n}$.
\end{definition}
The characteristic function can be also viewed as a homomorphism
$\l:\Z^m\to\Z^n$, where $m=\#{\mathcal F}=f_0$ and $\Z^m$ is the free
$\Z$-module spanned by the elements of ${\mathcal F}$.

It follows from the local equivalence of the torus action to the
standard one that the characteristic function has the following property: if
$F_{i_1},\ldots,F_{i_n}$ are the codimension-one faces meeting at the same
vertex, then $\l(F_{i_1}),\ldots,\l(F_{i_n})$ form an integer basis of
$\Z^n$. For any function $\l:{\mathcal F}\to\Z^n$ satisfying this condition
there exists a quasitoric manifold $M^{2n}(\l)$ over $P^n$ with
characteristic function $\l$, and $M^{2n}$ is determined by its
characteristic function up to an equivariant homeomorphism. However,
there are simple polytopes that do not admit any characteristic function.
One of such examples is the dual to the so-called cyclic polytope $C^n_k$
for $k\ge 2^n$ (see~\cite{DJ}). This polytope can not be realized as the
orbit space for a quasitoric (or toric) manifold.

Quasitoric manifolds over a simple polytope $P^n$ are
connected with the spaces $\ZP$ and $B_TP$ (see the previous
subsection) as follows. Since a quasitoric manifold $M^{2n}$ is a
$T^n$-space, we can take the Borel construction $ET^n\times_{T^n}M^{2n}$. It
turns out that all these spaces for the given combinatorial type $P^n$ are
independent of $M^{2n}$ and have the homotopy type $B_TP$:
\begin{equation}
\label{torbun}
  B_TP\simeq ET^n\times_{T^n}M^{2n}.
\end{equation} The relationship between the
$T^m$-space $\ZP$ and quasitoric manifolds over $P^n$ is described by the
following property: for each toric manifold $M^{2n}$ over $P^n$ the orbit map
$\ZP\to P^n$ is decomposed as $\ZP\to
M^{2n}\stackrel{\pi}{\longrightarrow}P^n$, where $\ZP\to M^{2n}$ is a
principal $T^{m-n}$-bundle and $M^{2n}\stackrel{\pi}{\longrightarrow}P^n$ is
the orbit map for $M^{2n}$. Therefore, quasitoric manifolds
over the polytope $P^n$ correspond to rank $m-n$ subgroups
of $T^m$ that act freely on $\ZP$; and each subgroup of this type produces
a quasitoric manifold. As it follows from Proposition~\ref{combstr},
subgroups of rank $\ge m-n$ can not act on $\ZP$ freely, so the
action of $T^m$ on $\ZP$ is maximally free  exactly when there is at least
one quasitoric manifold over $P^n$. We will discuss this
question in more details later (see subsection~\ref{addprop}).

From (\ref{bun}) we obtain the bundle $p:B_TP\to BT^m$ with the fibre $\ZP$.
All the cohomologies below are with the coefficients in the ground ring $k$.

\begin{theorem}
\label{sur}
  Let $P$ be a simple polyhedral complex with $m$ codimension-one faces.
  The map $p^*:H^*(BT^m)\to H^*(B_TP)$ is an epimorphism, and after the
  identification $H^*(BT^m)\cong k[v_1,\ldots,v_m]$ it becomes the
  quotient epimorphism $k[v_1,\ldots,v_m]\to k(P)$, where $k(P)$ is the face
  ring. In particular, $H^*(B_TP)\cong k(P)$. \ep
\end{theorem}

Now let $M^{2n}$ be a quasitoric manifold over a simple polytope $P^n$ with
the characteristic function $\l$. The characteristic function is obviously
extended to a linear map $k^m\to k^n$. Consider the bundle $p_0:B_TP\to
BT^n$ with the fibre $M^{2n}$ (see~(\ref{torbun})).

\begin{theorem}
\label{in}
  The map $p_0^*:H^*(BT^n)\to H^*(B_TP)$ is a monomorphism and
  $p_0^*:H^2(BT^n)\to H^2(B_TP)$ coincides with $\l^*:k^n\to k^m$.
  Furthermore,
  after the identification $H^*(BT^n)\cong  k[t_1,\ldots,t_n]$, the elements
  $\l_i=p^*_0(t_i)\in H^*(B_TP)\cong k(P)$ form a regular sequence of
  degree-two elements of $k(P)$.\ep
\end{theorem}

In particular, all the above constructions hold for (algebraic) toric
manifolds. Toric manifolds are defined by simple polytopes $P^n\subset\R^n$
whose vertices have integer coordinates. As it follows from the above
arguments, the value of the corresponding characteristic function on the
facet $F^{n-1}\in\mathcal F$ is its minimal integer normal (co)vector. All
characteristic functions corresponding to (algebraic) toric manifolds can be
obtained by this method.

\section{Geometrical and homotopical properties of $\ZP$ and
$B_TP$}

\subsection{A cubical subdivision of a simple polytope.}

In this subsection we suppose that $P^n$ is a simple $n$-dimensional
polytope. Let $I^q$ denote the standard unit cube in $\R^q$:
$$
  I^q=\{(y_1,\ldots,y_q)\in{\mathbb R}^q\,:\;0\le y_i\le 1,\,i=1,\ldots,q\}.
$$

\begin{definition}
\label{cub}
  A $q$-dimensional {\it cubical complex}
  is a topological space $X$ presented as the union of homeomorphic images of
  $I^q$ (called {\it cubes}) in such a way that the intersection of any two
  cubes is a face of each.
\end{definition}

\begin{theorem}
\label{Pcub}
  A simple polytope $P^n$ with $m=f_0$ facets is
  naturally a $n$-dimen\-sional cubical complex $\mathcal C$ with
  $r=f_{n-1}$ cubes $I^n_v$ corresponding to the vertices $v\in
  P^n$. Furthermore, there is an embedding $i_P$ of $\mathcal C$ into
  the boundary of the standard $m$-dimensional cube $I^m$ that
  takes cubes of $\mathcal C$ to $n$-faces of $I^m$.
\end{theorem}
\begin{proof}
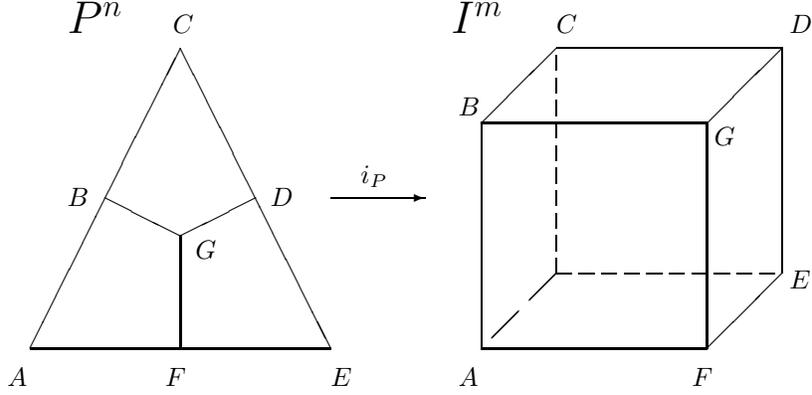
\begin{figure}
\begin{picture}(120,60)
  \put(10,10){\line(1,2){20}}
  \put(30,50){\line(1,-2){20}}
  \put(10,10){\line(1,0){40}}
  \put(20,30){\line(2,-1){10}}
  \put(40,30){\line(-2,-1){10}}
  \put(30,10){\line(0,1){15}}
  \put(50,30){\vector(1,0){12.5}}
  \put(70,10){\line(1,0){30}}
  \put(70,10){\line(0,1){30}}
  \put(70,40){\line(1,0){30}}
  \put(100,10){\line(0,1){30}}
  \put(70,40){\line(1,1){10}}
  \put(100,40){\line(1,1){10}}
  \put(100,10){\line(1,1){10}}
  \put(80,50){\line(1,0){30}}
  \put(110,20){\line(0,1){30}}
  \put(80,20){\line(-1,-1){3.6}}
  \put(75,15){\line(-1,-1){3.6}}
  \multiput(80,20)(3,0){10}{\line(1,0){2}}
  \multiput(80,20)(0,3){10}{\line(0,1){2}}
  \put(7,5){$A$}
  \put(28,5){$F$}
  \put(50,5){$E$}
  \put(32,22){$G$}
  \put(15,29){$B$}
  \put(42,29){$D$}
  \put(29,52){$C$}
  \put(15,52){\huge $P^n$}
  \put(54,32){$i_P$}
  \put(67,5){$A$}
  \put(98,5){$F$}
  \put(67,41){$B$}
  \put(101,37){$G$}
  \put(80,52){$C$}
  \put(111,52){$D$}
  \put(111,18){$E$}
  \put(66,52){\huge $I^m$}
\end{picture}
\caption{The embedding $i_p:P^n\to I^m$ for $n=2$, $m=3$.}
\label{fig1}
\end{figure}
Let us choose a point in the relative interior of each face of $P^n$ (we
also take all vertices and a point in the interior of the polytope). The
resulting set $\mathcal S$ of $1+f_0+f_1+\ldots+f_{n-1}$ points will be
the vertex set of the cubical complex $\mathcal C$. Since the polytope $P^n$
is simple, the number of $k$-faces meeting at each vertex is $\binom nk$,
$0\le k\le n$. Hence, for each vertex $v$ there is defined a $2^n$-element
subset $\mathcal S_v$ of $\mathcal S$ consisting of the points chosen in the
interiors of faces containing $v$ (including $v$ itself and the point in the
interior of $P^n$).
This set ${\mathcal S}_v$ is said to be the vertex set of the cube
$I^n_v\subset\mathcal C$ corresponding to $v$. The faces of $I^n_v$ are
defined as follows. We take any two faces $F_1^k$ and $F_2^l$ of $P^n$ such
that $v\in F_1^k\subset F_2^l$, $0\le k=\dim F^k\le l=\dim F^l\le n$. Then
there are $\binom{l-k}i$ faces $F^{k+i}$ of dimension $k+i$ such that $v\in
F_1^k\subset F^{k+i}\subset F_2^l$, $0\le i\le l-k$. Hence, there are
$2^{l-k}$ faces ``between" $F_1^k$ and $F_2^l$. The points inside these faces
define a $2^{l-k}$-element subset of $\mathcal S_v$, which is said to be a
vertex set of a $(l-k)$-face $I^{l-k}_{F_1,F_2}$ of the cube $I^n_v$. Now, to
finish the definition of the cubical complex $\mathcal C$ we need only to
check that the intersection of any two cubes $I^n_v$, $I^n_{v'}$ is a face of
each. To do this we take the minimal-dimension face $F^p\subset P^n$ that
contains both vertices $v$ and $v'$ (clearly, there is only one such face).
Then it can be easily seen that $I^n_v\cap I^n_{v'}=I^{n-p}_{F^p,P^n}$ is the
face of $I^n_v$ and $I^n_{v'}$.

Now let us construct an embedding ${\mathcal C}\hookrightarrow I^m$. First,
we define the images of the vertices of $\mathcal C$, i.e. the images of
the points of $\mathcal S$. To do this, we fix the numeration of facets:
$F^{n-1}_1,\ldots,F^{n-1}_m$. Now, if a point of $\mathcal S$ lies inside the
facet $F^{n-1}_i$, then we map it to the vertex $(1,\ldots,1,0,1,\ldots,1)$
of the cube $I^m$, where 0 stands on the $i$th place. If a point of
$\mathcal S$ lies inside a face $F^{n-k}$ of codimension $k$, then we write
$F^{n-k}=F^{n-1}_{i_1}\cap\ldots\cap F^{n-1}_{i_k}$, and map this point
to the vertex of $I^m$ whose $y_{i_1},\ldots,y_{i_k}$ coordinates are zero
and all other coordinates are 1. The point of $\mathcal S$ in the interior of
$P^n$ maps to the vertex of $I^m$ with coordinates $(1,\ldots,1)$. Hence, we
constructed the map from the set $\mathcal S$ to the vertex set of $I^m$.
This map obviously extends to a map from the cubical subdivision $\mathcal C$
of $P^n$ to the standard cubical subdivision of $I^m$. One of the ways to do
this is as follows. Take a
simplicial subdivision $\mathcal K$ of $P^n$ with vertex set $\mathcal S$
such that for each vertex $v\in P^n$ there exists a simplicial subcomplex
$\mathcal K_v\subset\mathcal K$ with vertex set $\mathcal S_v$ that
subdivides the cube $I^n_v$. The simplest way to construct such a simplicial
complex is to view $P^n$ as the cone over the barycentric subdivision of the
complex $K_P^{n-1}$ dual to the boundary $\partial P$.
Then the subcomplexes $\mathcal K_v$
are just the cones over the barycentric subdivisions of the $(n-1)$-simplices
of $K_P^{n-1}$.  Now we can extend the map ${\mathcal S}\hookrightarrow I^m$
linearly on each simplex of the triangulation $\mathcal K$ to the embedding
$i_P:P^n\hookrightarrow I^m$ (which is therefore a piecewise linear map).
Figure~\ref{fig1} illustrates this embedding for $n=2$, $m=3$.

In short, the above constructed embedding
$i_P:P^n\hookrightarrow I^m$ is determined by the following property:
$$
\begin{array}{c}
  \parbox{0.9\textwidth}
  {\noindent The cube $I^n_v\subset P^n$ corresponding to a vertex
  $v=F^{n-1}_{i_1}\cap\cdots\cap F^{n-1}_{i_n}$
  is mapped onto the $n$-face of $I^m$ determined
  by $m-n$ equations $y_j=1$, $j\notin \{i_1,\ldots,i_n\}$.}
\end{array}
\eqn
\label{cubmap}
$$
Thus, all cubes of $\mathcal C$ map to faces of $I^m$, which proves the
assertion.
\end{proof}

\begin{lemma}
  The number $c_k$ of $k$-cubes in the cubical subdivision $\mathcal C$
  of a simple polytope $P^n$ is given by the formula
  $$
    c_k=\sum_{i=0}^{n-k}f_{n-i-1}\binom{n-i}k=
    f_{n-1}\binom nk+f_{n-2}\binom{n-1}k+\ldots+f_{k-1},
  $$
  where $(f_0,\ldots,f_{n-1})$ is the $f$-vector of $P^n$ and $f_{-1}=1$.
\end{lemma}
\begin{proof}
This follows from the fact that $k$-cubes of
$\mathcal C$ are in one-to-one correspondence with pairs
$(F_1^{i},F_2^{i+k})$ of faces of $P^n$ such that $F_1^{i}\subset F_2^{i+k}$
(see the proof of Theorem~\ref{Pcub}).
\end{proof}

\subsection{$\ZP$ as a smooth manifold and an equivariant embedding of
$\ZP$ in $\C^m$.}

Let us consider the standard polydisc $(D^2)^m\subset\C^m$:
$$
  (D^2)^m=\{(z_1,\ldots,z_m)\in\C^m:\;\:|z_i|\le 1\}.
$$
This $(D^2)^m$ is a $T^m$-stable subset of $\C^m$ (with respect to the
standard action of $T^m$ on $\C^m$ by diagonal matrices). The corresponding
orbit space is $I^m$. The main result of this
subsection is the following theorem.

\begin{theorem}
\label{manif}
  For any simple polytope $P^n$ with $m$ facets the space $\ZP$ has the
  canonical structure of a smooth $(m+n)$-dimensional manifold for which the
  $T^{m}$-action is smooth. Furthermore, there exists a
  $T^m$-equivariant embedding $i_e:{\mathcal Z}_P\hookrightarrow
  (D^2)^m\subset\C^m$.
\end{theorem}
\begin{proof}
Theorem \ref{Pcub} shows that $P^n$ is presented as the union of $n$-cubes
$I^n_v$ indexed by the vertices of $P^n$.  Let $\rho:\ZP\to P^n$ be the orbit
map.  It follows easily from the definition of $\ZP$ that for each cube
$I^n_v\subset\ZP$ we have $\rho^{-1}(I^n_v)\cong(D^2)^n\times T^{m-n}$, where
$(D^2)^n$ is the polydisc in $\C^n$ with the diagonal action of $T^n$. Hence,
$\ZP$ is presented as the union of ``blocks" of the form $B_v\cong
(D^2)^n\times T^{m-n}$. Gluing these ``blocks" together smoothly along their
boundaries we obtain a smooth structure on $\ZP$. Since each $B_v$ is
$T^m$-invariant, the $T^m$-action on $\ZP$ is also smooth.

Now, let us prove the second part of the theorem.
Recall our numeration of codimension-one faces of $P^n$:
$F^{n-1}_1,\ldots,F^{n-1}_m$. Take the block
$$
  B_v\cong(D^2)^n\times T^{m-n}=
  D^2\times\ldots\times D^2\times S^1\times\ldots\times S^1
$$
corresponding to a vertex $v\in P^n$. Each factor $D^2$ and $S^1$ above
corresponds to a codimension-one face of $P^n$ and therefore acquires a
number (index) $i$, $1\le i\le m$. Note
that $n$ factors $D^2$ acquire the indices of those facets
containing $v$, while other indices are assigned to
$m-n$ factors $S^1$. Now we numerate the factors $D^2\subset(D^2)^m$ of the
polydisc in any way and embed each block $B_v\subset\ZP$ into $(D^2)^m$
according to the indices of its factors. It can be easily seen the set of
embeddings $B_v\hookrightarrow (D^2)^m$ define an equivariant embedding
$\ZP\hookrightarrow (D^2)^m$.
\end{proof}

\begin{example}
\label{sphdec}
  If $P^n=\Delta^1$ is a 1-dimensional simplex (a segment), then
  $B_v=D^2\times S^1$ for each of the two vertices, and we obtain the
  well-known decomposition
  $\mathcal Z_{\D^1}\cong S^3=D^2\times S^1\cup D^2\times
  S^1$. If $P^n=\Delta^n$ is a $n$-dimensional simplex, we obtain the
  similar decomposition of a $(2n+1)$-sphere into $n+1$ ``blocks"
  $(D^2)^n\times S^1$.
\end{example}

\begin{lemma}
\label{pullback}
  The equivariant embedding $i_e:\ZP\hookrightarrow (D^2)^m\subset\C^m$
  (see {\rm Theorem~\ref{manif}}) covers the embedding
  $i_P:P^n\hookrightarrow I^m$ (see {\rm Theorem~\ref{Pcub}}) as
  described by the commutative diagram
  $$
  \begin{CD}
    \ZP @>i_e>> (D^2)^m\\
    @VVV @VVV\\
    P^n @>i_P>> I^m,\\[2mm]
  \end{CD}
  $$
  where the vertical arrows denote the orbit maps for the corresponding
  $T^m$-actions.
\end{lemma}
\begin{proof}
It can be easily seen that an embedding of a face
$I^n\subset I^m$ defined by $m-n$ equations of the type $y_j=1$ (as
in~(\ref{cubmap})) induces an equivariant embedding
of $(D^2)^n\times T^{m-n}$ into $(D^2)^m$. Then our assertion follows from the
representation of $\ZP$ as the union of blocks
$B_v\cong(D^2)^n\times T^{m-n}$ and from property~(\ref{cubmap}).
\end{proof}

The above constructed embedding $i_e:\ZP\hookrightarrow (D^2)^m\subset\C^m$
allows as to relate the manifold $\ZP$ with one construction from the theory
of toric varieties. Below we describe this construction, following~\cite{Ba}.

We introduce the complex $m$-dimensional space $\C^m$ with coordinates
$z_1,\ldots,z_m$.

\begin{definition}
\label{primcol}
  A subset of facets $\mathcal
  P=\{F_{i_1},\ldots,F_{i_p}\}\subset\mathcal F$ is called a {\it primitive
  collection} if $F_{i_1}\cap\ldots\cap F_{i_p}=\emptyset$, while any
  proper subset of $\mathcal P$ has non-empty
  intersection. In terms of the simplicial complex $K_P$ dual to the
  boundary of $P^n$, the vertex subset
  $\mathcal P=\{v_{i_1},\ldots,v_{i_p}\}$ is called a primitive collection
  if $\{v_{i_1},\ldots,v_{i_p}\}$ does not span a simplex, while any proper
  subset of $\mathcal P$ spans a simplex of $K_P$.
\end{definition}

Now let $\mathcal P=\{F_{i_1},\ldots,F_{i_p}\}$ be a primitive collection
of facets of $P^n$. Denote by $\mathbf A(\mathcal P)$ the $(m-p)$-dimensional
affine subspace in $\mathbb C^m$ defined by the equations
\begin{equation}
\label{hplane}
  z_{i_1}=\ldots=z_{i_p}=0.
\end{equation}
Since every primitive collection has at least two elements, the
codimension of $\mathbf A(\mathcal P)$ is at least 2.

\begin{definition}
\label{defU}
  {\it The set of planes} (or {\it the arrangement of planes\/})
  $\mathbf A(P^n)\subset\mathbb C^m$ {\it defined by
  the lattice of faces of a simple polytope $P^n$\/} is
  $$
    \mathbf A(P^n)=\bigcup_{\mathcal P}\mathbf A(\mathcal P),
  $$
  where the union is taken over all primitive collections of facets of
  $P^n$. Put
  $$
    U(P^n)=\mathbb C^m\setminus \mathbf A(P^n).
  $$
\end{definition}

Note that we may define $U(P^n)$ without using primitive collections: the
same set is obtained if we take the complement in $\mathbb C^m$ of the union
of {\it all} planes~(\ref{hplane}) such that the facets
$F_{i_1},\ldots,F_{i_p}$ have empty intersection. However, we shall use
the notion of a primitive collection in the next sections. We note also that
the open set $U(P^n)\subset\mathbb C^m$ is invariant with respect to the
action of $(\mathbb C^{*})^m$ on $\mathbb C^m$.

It follows from property~(\ref{cubmap}) that the image of $\ZP$ under
the embedding $i_e:\ZP\to\C^m$ (see Theorem~\ref{manif}) does not intersect
$\mathbf A(P^n)$, and therefore, $i_e(\ZP)\subset U(P^n)$.

We put
$$
  \R^m_>=\{(\alpha_1,\ldots,\alpha_m)\in\R^n:\alpha_i>0\}.
$$
This is a group with respect to multiplication, which acts by dilations on
$\R^m$ and $\C^m$ (an element
$(\alpha_1,\ldots,\alpha_m)\in\R^m_>$ takes $(y_1,\ldots,y_m)\in\R^m$ to
$(\alpha_1y_1,\ldots,\alpha_my_m)$).
There is the isomorphism $\exp:\R^m\to\R^m_>$ between the additive and the
multiplicative group taking $(t_1,\ldots,t_m)\in\R^m$ to
$(e^{t_1},\ldots,e^{t_m})\in\R^m_>$.

Now, by definition, the polytope $P^n$ is a set of points $x\in\R^n$
satisfying $m$ linear inequalities:
\begin{equation}
\label{ptope}
  P^n=\{x\in\R^n:\<l_i,x\>\ge-a_i,\; i=1,\ldots,m\},
\end{equation}
where $l_i\in(\R^n)^*$ are normal (co)vectors of facets. The set
of $(\mu_1,\ldots,\mu_m)\in\R^m$ such that $\mu_1l_1+\ldots+\mu_ml_m=0$ forms
an $(m-n)$-dimensional subspace in $\R^m$. We choose a basis
$\{w_i=(w_{1i},\ldots,w_{mi})^\top\}$, $1\le i\le m-n$, in this subspace and
form the $m\times(m-n)$-matrix
\begin{equation}
\label{wmatrix}
  W=\begin{pmatrix}
  w_{11}&\ldots&w_{1,m-n}\\
  \ldots&\ldots&\ldots\\
  w_{m1}&\ldots&w_{m,m-n}
\end{pmatrix}
\end{equation}\\
of maximal rank $m-n$. This matrix satisfies the following property.

\begin{proposition}
\label{fprop}
  Suppose that $n$ facets $F^{n-1}_{i_1},\ldots,F^{n-1}_{i_n}$ of $P^n$
  meet at the same vertex $v$: $F^{n-1}_{i_1}\cap\cdots\cap F^{n-1}_{i_n}=v$.
  Then the minor $(m-n)\times(m-n)$-matrix $W_{i_1\ldots i_n}$ obtained
  from $W$ by deleting $n$ rows $i_1,\ldots,i_n$ is non-degenerate:
  $\det W_{i_1\ldots i_n}\ne0$.
\end{proposition}
\begin{proof}
If $\det W_{i_1,\ldots,i_n}=0$,
then one can find a zero non-trivial linear combination of
vectors $l_{i_1},\ldots,l_{i_n}$. But this is impossible: since $P^n$ is
simple, the set of normal vectors of facets meeting at the
same vertex constitute a basis of $\R^n$.
\end{proof}

The matrix $W$ defines the subgroup
$$
  R_W=\{(e^{w_{11}\tau_1+\cdots+w_{1,m-n}\tau_{m-n}},\ldots,
  e^{w_{m1}\tau_1+\cdots+w_{m,m-n}\tau_{m-n}})\in\R^m_>\}\subset\R^m_>,
$$
where $(\tau_1,\ldots,\tau_{m-n})$ runs over $\R^{m-n}$. This subgroup is
isomorphic to $\R^{m-n}_>$. Since $U(P^n)\subset\C^m$ (see
Definition~\ref{defU}) is invariant with respect to the
action of $\R^m_>\subset(\mathbb C^{*})^m$ on $\C^m$, the subgroup
$R_W\subset\R_>^m$ also acts on $U(P^n)$.

\begin{theorem}
\label{zu}
The subgroup $R_W\subset\R^m_>$ acts on $U(P^n)\subset\C^m$ freely.
The composite map $\ZP\to U(P^n)\to U(P^n)/R_W$ of the
embedding $i_e$ and the orbit map is a homeomorphism.
\end{theorem}
\begin{proof}
A point from $\C^m$ may have the non-trivial isotropy
subgroup with respect to the action of $\R^m_>$ on $\C^m$ only if at least
one its coordinate vanish. As it follows from Definition~\ref{defU},
if a point $x\in U(P^n)$ has some zero coordinates, then all of them
correspond to facets of $P^n$ having at least
one common vertex $v\in P^n$. Let $v=F^{n-1}_{i_1}\cap\cdots\cap
F^{n-1}_{i_n}$. The point $x$ has non-trivial isotropy subgroup with
respect to the action of $R_W$ only if some linear
combination of vectors $w_1,\ldots,w_{m-n}$ lies in the coordinate subspace
spanned by $e_{i_1},\ldots,e_{i_n}$. But this means that
$\det W_{i_1\ldots i_n}=0$, which contradicts
Proposition~\ref{fprop}. Thus, $R_W$ acts on $U(P^n)$ freely.

Now, let us prove the second part of the theorem. Here we use both
embeddings $i_e:\ZP\to(D^2)^m\subset\C^m$ from Theorem~\ref{manif} and
$i_P:P^n\to I^m\subset\R^m$ from Theorem~\ref{Pcub}. It is sufficient to
prove that each orbit of the action of $R_W$ on $U(P^n)\subset\C^m$
intersects the image $i_e(\ZP)$ in a single point. Since the embedding $i_e$
is equivariant, instead of this we may prove that each orbit of the action of
$R_W$ on the real part $U_{\R}(P^n)=U(P^n)\cap\R^m_+$ intersects the image
$i_P(P^n)$ in a single point. Let $y\in i_P(P^n)\subset\R^m$. Then
$y=(y_1,\ldots,y_m)$ lies in some $n$-face $I^n_v$ of the unit cube
$I^m\subset\R^m$ as described in~(\ref{cubmap}). We need to show that the
$(m-n)$-dimensional subspace spanned by the vectors
$(w_{11}y_1,\ldots,w_{m1}y_m)^\top,\ldots,
(w_{1,m-n}y_1,\ldots,w_{m,m-n}y_m)^\top$ is in general position with the
$n$-face $I^n_v$. But this follows directly from~(\ref{cubmap}) and
Proposition~\ref{fprop}.
\end{proof}

The above theorem gives a new proof of the fact that $\ZP$ is a smooth
manifold, which is embedded in $\C^m\cong\R^{2m}$ with trivial normal bundle.

\begin{example}
\label{sphere}
  Let $P^n=\D^n$ ($n$-simplex). Then
  $m=n+1$, $U(P^n)=\C^{n+1}\setminus\{0\}$, $R^{m-n}_>$ is
  $\R_>$, and $\alpha\in\R_>$ takes $z\in\C^{n+1}$ to $\alpha z$.
  Thus, we have $\ZP=S^{2n+1}$ (this could be also deduced
  from Definition~\ref{defzp}; see also Example~\ref{sphdec}).
\end{example}

Now, suppose that all vertices of $P^n$ belong to the integer lattice
$\Z^n\subset\R^n$. Such integer simple polytope $P^n$ defines a
projective toric variety $M_P$~(see~\cite{Fu}). Normal (co)vectors $l_i$
of facets of $P^n$ (see~(\ref{ptope})) can be taken integer and primitive.
The toric variety $M_P$ defined by $P^n$ is smooth if
for each vertex $v=F_{i_1}\cap\ldots\cap F_{i_n}$ the vectors
$l_{i_1},\ldots,l_{i_n}$ constitute an integer basis of $\Z^n$.
As before, we may construct the matrix $W$ (see~(\ref{wmatrix})) and
then define the subgroup
$$
  C_W=\{(e^{w_{11}\tau_1+\cdots+w_{1,m-n}\tau_{m-n}},\ldots,
  e^{w_{m1}\tau_1+\cdots+w_{m,m-n}\tau_{m-n}})\}\subset
  (\C^{*})^m,
$$
where $(\tau_1,\ldots,\tau_{m-n})$ runs over $\C^{m-n}$. This subgroup is
isomorphic to $(\C^{*})^{m-n}$. It can be shown (see~\cite{Ba}) that $C_W$
acts freely on $U(P^n)$ and the toric manifold $M_P$ is identified with
the orbit space $U(P^n)/C_W$. Thus, we have the commutative diagram
$$
\begin{CD}
  U(P^n) @>R_W\cong\R^{m-n}_{>}>> \ZP\\
  @VC_W\cong(\C^{*})^{m-n}VV @VVT^{m-n}V\\
  M^{2n} @= M^{2n}.\\[2mm]
\end{CD}
$$

Since $\ZP$ can be viewed as the orbit space of $U(P^n)$ with
respect to an action of $R_W\cong\R^{m-n}_>$, the manifold $\ZP$ and the
complement of an arrangement of planes $U(P^n)$ have the same homotopy type.
Hence, all results on the cohomology of $\ZP$ obtained in section~4
remain true if we substitute $U(P^n)$ for $\ZP$.

\subsection{Homotopical properties of $\ZP$ and $B_TP$.}

We start with two simple assertions.

\begin{lemma}
\label{prod}
  Suppose that $P^n$ is the product of two simple polytopes,
  $P^n=P^{n_1}_1\times P^{n_2}_2$. Then
  $\ZP={\mathcal Z}_{P_1}\times{\mathcal Z}_{P_2}$.
\end{lemma}
\begin{proof}
This follows directly from the definition of $\ZP$:
$$
  \ZP=(T^m\times P^n)/\!\sim\:=
  \bigl((T^{m_1}\times P^{n_1})/\!\sim\bigr)
  \times\bigl((T^{m_2}\times P^{n_2})/\!\sim\bigr)=
  {\mathcal Z}_{P_1}\times{\mathcal Z}_{P_2}.
$$
\end{proof}
The next lemma also follows easily from the construction of $\ZP$.
\begin{lemma}
\label{face}
  If $P_1^{n_1}\subset P^n$ is a face of a simple polytope
  $P^n$, then ${\mathcal Z}_{P_1}$ is a submanifold of $\ZP$.\ep
\end{lemma}

Below we invest the space $B_TP$ defined by a simple polyhedral complex
with a canonical cell structure.

We use the standard cell decomposition of $BT^m=(\C P^{\infty})^m$ (each
$\C P^{\infty}$ has one cell in every even dimension).
The corresponding cell cochain algebra is
$C^{*}(BT^m)=H^{*}(BT^m)=k[v_1,\ldots,v_m]$.

\begin{theorem}
\label{cell}
  The space $B_TP=ET^m\times_{T^m}\ZP$ defined by a simple polyhedral complex
  $P$ can be viewed as a cell subcomplex of $BT^m$. This subcomplex
  is the union of subcomplexes $BT^k_{i_1,\ldots,i_k}$ over all
  simplices $\Delta=(i_1,\ldots,i_k)$ of the simplicial complex $K_P^{n-1}$
  dual to the boundary $\partial P^n$. In this realization we have
  $C^{*}(B_TP)=H^{*}(B_TP)=k(P)$, and the inclusion
  $i:B_TP\hookrightarrow BT^m$ induces the quotient epimorphism
  $C^{*}(BT^m)=k[v_1,\ldots,v_m]\to k(P)=C^{*}(B_TP)$ (here
  $k(P)$ is the face ring of $P$).
\end{theorem}
\begin{proof}
A simple polyhedral complex $P$ is defined
as the cone over the barycentric subdivision of a
simplicial complex $K$ with $m$ vertices. We construct a cell embedding
$i:B_TP\hookrightarrow BT^m$ by induction on the dimension of $K$.
If $\dim K=0$, then $K$ is a disjoint union of vertices $v_1,\ldots,v_m$ and
$P$ is the cone on $K$. In this case $B_TP$ is a bouquet of $m$ copies of
$\C P^{\infty}$ and we have the obvious inclusion
$i:B_TP\to BT^m=(\C P^{\infty})^m$. In degree zero $C^{*}(B_TP)$ is just $k$,
while in degrees $\ge 1$ it is isomorphic to
$k[v_1]\oplus\cdots\oplus k[v_m]$. Therefore,
$C^{*}(B_TP)=k[v_1,\ldots,v_m]/I$, where $I$ is the ideal generated by all
square free monomials of degree $\ge 2$, and $i^{*}$ is the projection onto
the quotient ring. Thus, the theorem holds for $\dim K=0$.

Now let $\dim K=k-1$. By the inductive hypothesis, the theorem holds for the
simple polyhedral complex $P'$ corresponding to the $(k-2)$-skeleton $K'$ of
$K$, i.e.  $i^{*}C^{*}(BT^m)=C^{*}(B_TP')=k(K')=k[v_1,\ldots,v_m]/I'$. We add
${(k-1)}$-simplices one at a time. Adding the simplex $\Delta^{k-1}$ on
vertices $v_{i_1},\ldots,v_{i_k}$ results in adding all cells of the
subcomplex $BT^k_{i_1,\ldots,i_k}=BT^1_{i_1}\times\ldots\times
BT^1_{i_k}\subset BT^m$ to $B_TP'\subset BT^m$. Then $C^{*}(B_TP'\cup
BT^k_{i_1,\ldots,i_k})= k(K'\cup\Delta^{k-1})=k[v_1,\ldots,v_m]/I$, where $I$
is generated by $I'$ and $v_{i_1}v_{i_2}\ldots v_{i_k}$. It is also clear
that a map of the cochain (or cohomology) algebras induced by
$i:B_TP\hookrightarrow BT^m$ is the projection onto the quotient ring.
\end{proof}

In particular, we see that for $K_P=\Delta^{m-1}$ (i.e. $P=\R^m_{+}$) one
has $B_TP=BT^m$.

Below we apply the above constructed cell decomposition of $B_TP$ for
calculating some homotopy groups of $B_TP$ and $\ZP$.

A simple polytope (or a simple polyhedral complex) $P^n$ with $m$
codimension-one faces is called {\it $q$-neighbourly}~\cite{Br} if the
$(q-1)$-skeleton of the dual simplicial complex $K^{n-1}_P$
coincides with the $(q-1)$-skeleton of a $(m-1)$-simplex
(this just means that any $q$ codimension-one faces of $P^n$ have non-empty
intersection). Note that any simple polytope is 1-neighbourly.

\begin{theorem}
\label{homot}
  For any simple polyhedral complex $P^n$ with $m$ codimension-one faces we
  have:
  \begin{enumerate}

  \item $\pi_1(\ZP)=\pi_1(B_TP)=0$;

  \item $\pi_2(\ZP)=0,\;\pi_2(B_TP)=\Z^m$;

  \item $\pi_q(\ZP)=\pi_q(B_TP)$ for $q\ge 3$;

  \item If $P^n$ is $q$-neighbourly, then
  $\pi_i(\ZP)=0$ for $i<2q+1$,
  and $\pi_{2q+1}(\ZP)$ is a free Abelian group with generators
  corresponding to square-free
  monomials $v_{i_1}\cdots v_{i_{q+1}}\in I$ (see
  {\rm Definition~\ref{frpol};} these monomials correspond to primitive
  collections of $q+1$ facets).
  \end{enumerate}
\end{theorem}
\begin{proof}
The identities $\pi_1(B_TP)=0$ and $\pi_2(B_TP)=\Z^m$ follow
from the cell decomposition of $B_TP$ described in the previous theorem. In
order to calculate $\pi_1({\mathcal Z}_P)$ and $\pi_2({\mathcal Z}_P)$
we consider the following fragment of the exact homotopy sequence
for the bundle $p:B_TP\to BT^m$ with the fibre $\ZP$:
$$
\begin{array}{lcccr}
  \pi_3(BT^m)\to\pi_2({\mathcal Z}_P)\to&\pi_2(B_TP)&
  \stackrel{p_{*}}{\longrightarrow}&\pi_2(BT^m)&\to\pi_1({\mathcal Z}_P)
  \to\pi_1(B_TP)\\[1mm]
  \quad\|&\|\;&&\;\|&\|\quad\\[1mm]
  \quad0&\Z^m&\longrightarrow&\Z^m&0\quad
\end{array}
$$
It follows from Theorem \ref{cell} that $p_{*}$ above is an isomorphism,
and hence, $\pi_1(\ZP)=\pi_2(\Z_P)=0$. The third assertion of the
theorem follows from the fragment
$$
  \pi_{q+1}(BT^m)\to\pi_q({\mathcal Z}_P)\to\pi_q(B_TP)\to\pi_q(BT^m),
$$
in which $\pi_q(BT^m)=\pi_{q+1}(BT^m)=0$ for $q\ge3$. Finally, the cell
structure of $B_TP$ shows that if $P^n$ is $q$-neighbourly, then
the $(2q+1)$-skeleton of $B_TP$ coincides with the $(2q+1)$-skeleton
of $BT^m$. Thus, $\pi_k(B_TP)=\pi_k(BT^m)$ for $k<2q+1$. Now, the last
assertion of the theorem follows from the third one and
from Theorem~\ref{cell}.
\end{proof}

The form of the homotopy groups of $\ZP$ and $B_TP$
enables us to make a hypothesis that
$\ZP$ is a first killing space for $B_TP$, i.e. $\ZP=B_TP|_{3}$.
In order to see this, let us
consider the following commutative diagram of bundles:
\begin{equation}
\begin{CD}
  \ZP\times ET^m @>>> ET^m\\
  @VVV @VVV\\
  B_TP @>p>> BT^m.
\end{CD}
\label{kill}
\end{equation}
Since $ET^m$ is contractible, $\ZP\times ET^m$ is homotopically equivalent to
$\ZP$. On the other hand, since $BT^m=K(\Z^m,2)$ and $\pi_2(B_TP)=\Z^m$,
we see that $\ZP\times ET^m$ is a first killing space for $B_TP$ by
definition. Thus, $\ZP$ has homotopy type of a first killing space
$B_TP|_{3}$ for $B_TP$.

\section{The Eilenberg--Moore spectral sequence.}

In \cite{EM} Eilenberg and Moore developed a spectral sequence, which
turns out to be of great use in our considerations. In
the description of this spectral sequence we follow~\cite{Sm}.

Suppose that $\xi_0=(E_0,p_0,B_0,F)$ is a Serre fibre bundle, $B_0$ is simply
connected, and $f:B\to B_0$ is a continuous map. Then we can form the diagram
\begin{equation}
\begin{CD}
  F @= F\\
  @VVV @VVV\\
  E @>>> E_0\\
  @VpVV @VVp_0V\\
  B @>f>> B_0,
\end{CD}
\label{comsq}
\end{equation}
where $\xi=(E,p,B,F)$ is the induced fibre bundle. Under these assumptions the
following theorem holds

\begin{theorem}[\rm Eilenberg--Moore]
  There exists a spectral sequence of commutative algebras $\{E_r,d_r\}$ with
  \begin{enumerate}
    \item $E_r\Rightarrow H^*(E)$ (the spectral sequence converges to
    the cohomology of $E$),
    \item $E_2=\Tor_{H^*(B_0)}\bigl(H^*(B),H^*(E_0)\bigr)$. \ep
  \end{enumerate}
\end{theorem}
The Eilenberg--Moore spectral sequence lives in the second quadrant and the
differential $d_r$ has bidegree $(r,1-r)$. In the special case when $B=*$
is a point (hence, $E=F$ is the fibre of $\xi$) we have

\begin{corollary}
\label{onefib}
  Let $F\hookrightarrow E\to B$ be a fibration over the simply connected
  space $B$. There exists a spectral sequence of commutative algebras
  $\{E_r,d_r\}$ with
  \begin{enumerate}
    \item $E_r\Rightarrow H^*(E)$,
    \item $E_2=\Tor_{H^*(B)}\bigl(H^*(E),k\bigr)$. \ep
  \end{enumerate}
\end{corollary}

As the first application of the Eilenberg--Moore spectral sequence we
calculate the cohomology ring of a quasitoric manifold $M^{2n}$ over a simple
polytope $P^n$ (this was already done in~\cite{DJ} by means of other
methods). Along with the ideal $I$ (see Definition~\ref{frpol}) we define an
ideal $J\subset k(P)$ as $J=(\l_1,\ldots,\l_n)$, where $\l_i$ are the
elements of the face ring $k(P)$ defined by the characteristic function $\l$
of the manifold $M^{2n}$ (see Theorem~\ref{in}). As it follows from
Theorem~\ref{in}, $\l_i=\l_{i1}v_1+\l_{i2}v_2+\ldots+\l_{im}v_m$ are
algebraically independent elements of degree $2$ in $k(P)$, and $k(P)$ is a
finite-dimensional free $k[\l_1,\ldots,\l_n]$-module. The inverse image of
the ideal $J$ under the projection $k[v_1,\ldots,v_m]\to k(P)$ is the
ideal generated by $\l_i=\l_{i1}v_1+\ldots+\l_{im}v_m$ regarded as
elements of $k[v_1,\ldots,v_m]$. This inverse image will be also denoted
by $J$.

\begin{theorem}
\label{cohomM}
  The following isomorphism of rings holds
  for any quasitoric manifold $M^{2n}$:
  $$
    H^*(M^{2n})\cong k(P)/J=k[v_1,\ldots,v_m]\,/\,I{+}J.
  $$
\end{theorem}
\begin{proof}
Consider the Eilenberg--Moore spectral sequence of the fibration
$$
\begin{CD}
  M^{2n} @>>> B_TP\\
  @VVV @VV\mbox{\small$p_0$}V\\
  {*} @>>> BT^n
\end{CD}
$$
Theorem~\ref{in} gives the monomorphism
\begin{eqnarray*}
  H^*(BT^n)=k[t_1,\ldots,t_n]&\stackrel{p_0^*}{\longrightarrow}&
  H^*(B_TP)=k(P),\\[1mm]
  t_i&\longrightarrow&\l_i,
\end{eqnarray*}
such that $\mathop{\rm Im}p_0^*=k[\l_1,\ldots,\l_n]\subset k(P)$. The
$E_2$ term of the Eilenberg--Moore spectral sequence is
$$
  E_2^{*,*}=\Tor^{*,*}_{H^*(BT^n)}\bigl(H^*(B_TP),k\bigr)=
  \Tor^{*,*}_{k[\l_1,\ldots,\l_n]}\bigl(k(P),k\bigr).
$$
The right-hand side above is a bigraded $k$-module
(see~\cite{Ma},~\cite{Sm}).
The first (``external") grading arises from a projective
resolution of  $H^*(B_TP)$ as a $H^*(BT^n)$-module used in the definition of
the functor $\Tor$. The second (``internal") grading arises from the
gradings of $H^*(BT^n)$-modules which enter the resolution; we assume that
non-zero elements appear only in even internal degrees
(remember that $\deg\l_i=2$).
Since $k(P)$ is a free $k[\l_1,\ldots,\l_n]$-module, we have
$$
  \Tor^{*,*}_{k[\l_1,\ldots,\l_n]}\bigl(k(P),k\bigr)=
  \Tor^{0,*}_{k[\l_1,\ldots,\l_n]}\bigl(k(P),k\bigr)=
  k(P)\otimes_{k[\l_1,\ldots,\l_n]}k=k(P)/J.
$$
Therefore, $E_2^{0,*}=k(P)/J$ and $E_2^{-p,*}=0$ for $p>0$. Thus,
$E_2=E_{\infty}$ and $H^*(M^{2n})=k(P)/J$.
\end{proof}

As we have already mentioned, this theorem generalizes the well-known
Danilov--Jurkiewicz theorem for the cohomology ring of a non-singular
projective toric variety.

\begin{corollary}
\label{Mtor}
  $H^*(M^{2n})=\Tor_{k[\l_1,\ldots,\l_n]}\bigl(k(P),k\bigr)$. \ep
\end{corollary}

\section{Calculation of the cohomology of $\ZP$}

In this section we use the Eilenberg--Moore spectral sequence for describing
the cohomology ring of $\ZP$ in terms of the face ring $k(P)$. We also
obtain some additional results about this cohomology in the case when at
least one quasitoric manifold exists over the polytope $P$. Throughout this
section we assume that $k$ is a field.

\subsection{Additive structure of the cohomology of $\ZP$.}

Here we consider the Eilenberg--Moore spectral
sequence of the bundle $p:B_TP\to BT^m$ with the fibre $\ZP$ (see~(\ref{bun})).
This spectral sequence defines a decreasing filtration on $H^{*}(\ZP)$, which
we denote $\{F^{-p}H^{*}({\ZP})\}$, such that
$$
  E_{\infty}^{-p,n+p}=F^{-p}H^{n}(\ZP)/F^{-p+1}H^{n}(\ZP).
$$

\begin{proposition}
\label{f0}
  $F^0H^{*}({\ZP})=H^0({\ZP})=k$ (here $k$ is the ground field).
\end{proposition}
\begin{proof}
It follows from \cite[Proposition 4.2]{Sm} that for the
Eilenberg--Moore spectral sequence of an arbitrary commutative
square~(\ref{comsq}) one has $F^0H^{*}(E)=\mathop{\rm Im}\{H^{*}(B)\otimes
H^{*}(E_0)\to H^{*}(E)\}$. In our case this gives
$F^0H^{*}({\ZP})=\mathop{\rm Im}\{H^{*}(B_TP)\to H^{*}({\ZP})\}$. Now, the
proposition follows from the fact that
the map $p^*:H^*(BT^m)\to H^*(B_TP)$ is an epimorphism (see
Theorem~\ref{sur}).
\end{proof}

\medskip

The $E_2$ term of the
Eilenberg--Moore spectral sequence of the bundle $p:B_TP\to BT^m$
is $E_2=\Tor_{k[v_1,\ldots,v_m]}\bigl(k(P),k\bigr)$.
Let us consider a free resolution of $k(P)$ as a $k[v_1,\ldots,v_m]$-module:
\begin{equation}
\label{resol}
  0\longrightarrow R^{-h}\stackrel{d^{-h}}{\longrightarrow}
  R^{-h+1}\stackrel{d^{-h+1}}{\longrightarrow}\cdots\longrightarrow
  R^{-1}\stackrel{d^{-1}}{\longrightarrow}
  R^{0}\stackrel{d^{0}}{\longrightarrow} k(P)\longrightarrow 0.
\end{equation}
It is convenient for our purposes to assume that $R^i$ are numbered by
non-positive integers, i.e. $h>0$ above.

The minimal number $h$ for which a free
resolution of the form~(\ref{resol}) exists is called the {\it homological
dimension} of $k(P)$ and is denoted by
$\mathop{\rm hd}\nolimits_{k[v_1,\ldots,v_m]}\bigl(k(P)\bigr)$.
By the Hilbert syzygy theorem,
$\mathop{\rm hd}\nolimits_{k[v_1,\ldots,v_m]}\bigl(k(P)\bigr)\le m$.
At the same time, since $k(P)$ is
a Cohen--Macaulay ring, it is known~\cite[Chapter IV]{Se} that
$$
  \mathop{\rm hd}\nolimits_{k[v_1,\ldots,v_m]}\bigl(k(P)\bigr)=m-n,
$$
where $n$ is the Krull dimension (the maximal number of algebraically
independent elements) of $k(P)$. In our case $n=\dim P$.

We shall use a special free resolution (\ref{resol}) known as
the {\it minimal} resolution (see~\cite{Ad}), which is defined in the
following way. Let $A$ be a graded connected commutative algebra, and let
$R,R'$ be modules over $A$. Set $I(A)=\sum_{q>0}A_q=\{a\in A:\:\deg a\ne 0\}$
and $J(R)=I(A)\cdot R$. The map $f:R\to R'$ is called {\it minimal}, if
${\mathop{\rm Ker}f\subset J(R)}$. The resolution (\ref{resol}) is called
{\it minimal}, if all $d^i$ are minimal. For constructing a minimal
resolution we use a notion of a {\it minimal set of generators}.
A minimal set of generators for a $A$-module $R$ can be chosen
by means of the following procedure.
Let $k_1$ be the lowest degree in which $R$ is non-zero.
Choose a vector space basis in $(R)^{k_1}$, say $x_1,\ldots,x_p$. Now let
$R_1=(x_1,\ldots,x_p)\subset R$ be the submodule generated by
$x_1,\ldots,x_p$. If $R=R_1$ then we are done.
Otherwise, consider the first degree $k_2$ in which $R\ne
R_1$; then in this degree we can choose a direct sum decomposition
$R=R_1\oplus\widehat{R_1}$. Now choose in $\widehat{R_1}$ a vector space
basis $x_{p_1+1},\ldots,x_{p_2}$ and set $R_2=(x_1,\ldots,x_{p_2})$. If
$R=R_2$ we are done, if not just continue to repeat the above process until
we obtain a minimal set of generators for $R$. A minimal set of generators
has the following property: no element $x_k$ can
be decomposed as $x_k=\sum a_ix_i$ with $a_i\in A$, $\deg a_i\ne 0$.
Now, we construct a minimal resolution~(\ref{resol}) as follows. Take a
minimal set of generators for $k(P)$ and span by them a free
$k[v_1,\ldots,v_m]$-module $R^0$. Then take a minimal set of generators for
$\mathop{\rm Ker}d^0$ and span by them a free module $R^{-1}$, and so on.
On the $i$th step we take a minimal set of homogeneous generators for
$\mathop{\rm Ker}d^{-i+1}$ as the basis for $R^{-i}$. A~minimal resolution is
unique up to an isomorphism.

Now let (\ref{resol}) be a minimal resolution of $k(P)$ as a
$k[v_1,\ldots,v_m]$-module. Then we have $h=m-n$, and $R^0$ is a free
$k[v_1,\ldots,v_m]$-module with one generator of degree $0$.  The generator
set of $R^1$ consists of elements $v_{i_1\ldots i_k}$ of degree $2k$ such
that $\{v_{i_1},\ldots,v_{i_k}\}$ does not span a simplex in $K$, while any
proper subset of $\{v_{i_1},\ldots,v_{i_k}\}$ spans a simplex in $K$. This
means exactly that the set $\{v_{i_1},\ldots,v_{i_k}\}$ is a primitive
collection in the sense of Definition~\ref{primcol} (here the $v_i$ are
regarded as the vertices of the simplicial complex $K^{n-1}$ dual to
$\partial P$).

Note that the $k[v_1,\ldots,v_m]$-module structure in $k$ is defined by the
homomorphism $k[v_1,\ldots,v_m]\to k$, $v_i\to 0$. Since the
resolution~(\ref{resol}) is minimal, all
the differentials $d^i$ in the complex
\begin{multline}
\label{timesk}
  \begin{CD}
  0 @>>> R^{-(m-n)}\otimes_{k[v_1,\ldots,v_m]}k
  @>d^{-(m-n)}>> \cdots
  \end{CD}\\
  \begin{CD}
  \cdots @>>> R^{-1}\otimes_{k[v_1,\ldots,v_m]}k
  @>d^{-1}>> R^{0}\otimes_{k[v_1,\ldots,v_m]}k @>>> 0
  \end{CD}
\end{multline}
are trivial. The module $R^i\otimes_{k[v_1,\ldots,v_m]}k$ is a
finite-dimensional vector space over $k$; its dimension is equal to the
dimension of $R^i$ as a free $k[v_1,\ldots,v_m]$-module:
$$
  \dim_kR^i\otimes_{k[v_1,\ldots,v_m]}k=\dim_{k[v_1,\ldots,v_m]}R^i.
$$
Therefore, since all the
differentials in the complex (\ref{timesk}) are trivial,
the following equality holds for the minimal resolution (\ref{resol}):
\begin{equation}
\label{dimtor}
  \dim_k\Tor_{k[v_1,\ldots,v_m]}\bigl(k(P),k\bigr)=
  \sum_{i=0}^{m-n}\dim_{k[v_1,\ldots,v_m]}R^{-i}.
\end{equation}

Now we are ready to describe the additive structure of the cohomology of
$\ZP$.

\begin{theorem}
\label{cohomZ}
  The following isomorphism of graded $k$-modules holds:
  $$
    H^*(\ZP)\cong\Tor_{k[v_1,\ldots,v_m]}\bigl(k(P),k\bigr),
  $$
  where the right-hand side is regarded as the
  totalized one-graded module.
  More precisely, there is a filtration $\{F^{-p}H^{*}(\ZP)\}$ in
  $H^*(\ZP)$ such that
  $$
    F^{-p}H^{*}(\ZP)/F^{-p+1}H^{*}(\ZP)=
    \Tor^{-p}_{k[v_1,\ldots,v_m]}\bigl(k(P),k\bigr).
  $$
\end{theorem}
\begin{proof}
First, we show that
\begin{equation}
\label{uneq}
  \dim_kH^{*}(\ZP)\ge\dim_k\Tor_{k[v_1,\ldots,v_m]}\bigl(k(P),k\bigr).
\end{equation}
To prove this inequality we construct an injective map from the union of
generator sets of the free $k[v_1,\ldots,v_m]$-modules $R^i$
(see~(\ref{resol})) to a generator set of $H^{*}(\ZP)$ (over $k$).

Consider the Leray--Serre spectral sequence of the bundle
$p:B_TP\to BT^m$ with the fibre $\ZP$. The first column of the $E_2$ term
of this spectral sequence is the cohomology of the fibre:
$H^{*}(\ZP)=E_2^{0,*}$.
By Theorem \ref{sur}, non-zero elements can appear in the $E_{\infty}$ term
only in the bottom line; this bottom line is the ring $k(P)=H^{*}(B_TP)$:
$$
  E_{\infty}^{*,p}=0,\;p>0;\quad E_{\infty}^{*,0}=k(P).
$$
Therefore, all elements from the kernel of the map
$d^0:\:R^0=k[v_1,\ldots,v_m]=E_2^{*,0}\to E_{\infty}^{*,0}=k(P)$
(see~(\ref{resol})) must be killed by the differentials of the
spectral sequence. This kernel is just our ideal $I$.
Let $(x_1,\ldots,x_p)$ be a minimal generator set of $I$. Then we claim that
the elements $x_i$ can be killed only by the transgression (i.e. by
differentials from the first column). Indeed, suppose that
the converse is true, so that
$x\in\{x_1,\ldots,x_p\}$ is killed by a
non-transgressive differential: $x=d_ky$ for some $k$, where $y$ is not from
the first column. Then $y$ is sent to zero by all differentials up to
$d_{k-1}$. This $y$ arises from some element $\sum_il_ia_i$ in the $E_2$
term, $l_i\in E_2^{0,*}$, $a_i\in E_2^{*,0}=k[v_1,\ldots,v_m]$. Suppose that
all the elements $l_i$ are transgressive (i.e. $d_j(l_i)=0$
for $j<k$), and let $d_k(l_i)=m_i$, $m_i\in E_k^{*,0}$.
Since all $m_i$ are killed by differentials, their inverse images in $E_2$
belong to $I$. Hence, we have $x=d_ky=\sum_im_ia_i$, $m_i\in I$, which
contradicts the minimality of the basis $(x_1,\ldots,x_p)$. Therefore,
there are some non-transgressive elements
among $l_i$, i.e. there exists $p<k$ and $i$ such that
$d_p(l_i)=m_i\ne0$ (see Figure~{\ref{fig2}}).
Then $m_i$ survives in $E_p$ and
we have $d_p(y)=m_ia_i+\ldots\ne0$~--- contradiction. This means that
all the minimal generators of $I$ are killed by the
transgression, i.e. they correspond to some (different) generators
$l_i^{(1)}\in H^{*}({\mathcal Z})$.

\begin{figure}
\begin{center}
\begin{picture}(60,40)
  \put(0,0){\line(0,1){30}}
  \put(0,0){\line(1,0){60}}
  \put(0,5){\line(1,0){60}}
  \put(5,0){\line(0,1){30}}
  \put(25,0){\line(0,1){30}}
  \put(30,0){\line(0,1){30}}
  \put(50,0){\line(0,1){5}}
  \put(55,0){\line(0,1){5}}
  \put(0,20){\line(1,0){60}}
  \put(0,25){\line(1,0){60}}
  \put(15,10){\line(1,0){5}}
  \put(15,10){\line(0,1){5}}
  \put(20,15){\line(-1,0){5}}
  \put(20,15){\line(0,-1){5}}
  \put(5,20){\vector(2,-1){10}}
  \put(30,20){\vector(4,-3){20}}
  \put(11.5,17){\footnotesize $d_p$}
  \put(41,12){\footnotesize $d_k$}
  \put(2,21.5){$l_i$}
  \put(15.3,12){$m_i$}
  \put(26,1.5){$a_i$}
  \put(26.5,22){$y$}
  \put(51.5,1.5){$x$}
\end{picture}%
\caption{The Leray--Serre spectral sequence of $p:B_TP\to BT^m$.}
\label{fig2}
\end{center}
\end{figure}
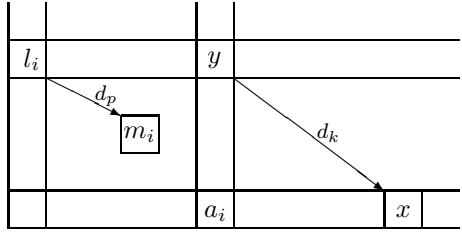

Since $E_2=H^{*}(\ZP)\otimes k[v_1,\ldots,v_m]$, a free
$k[v_1,\ldots,v_m]$-module generated by the elements $l_i^{(1)}$ is included
into the $E_2$ term as a submodule. Therefore, we have $R^{-1}\subset E_2$
and the map $d^{-1}:R^{-1}\to R^0=k[v_1,\ldots,v_m]$ is defined by the
differentials of the spectral sequence. The kernel of this map,
$\mathop{\rm Ker}d^{-1}$, can not be killed by the already constructed
differentials. Using the previous argument, we deduce that the elements of
a minimal generator set for $\mathop{\rm Ker}d^{-1}\in R^{-1}$ can be killed
only by some elements from the first column, say
$l_1^{(2)},\ldots,l_q^{(2)}$.  Therefore, a free $k[v_1,\ldots,v_m]$-module
generated by the elements $l_i^{(2)}$ is also included into the $E_2$ term as
a submodule, i.e.  $R^{-2}\subset E_2$. Proceeding with this procedure, at
the end we obtain $\sum_{i=0}^{m-n}\dim_{k[v_1,\ldots,v_m]}R^{-i}$ generators
in the first column of the $E_2$ term. Using~(\ref{dimtor}), we deduce the
required inequality~(\ref{uneq}).

Now let us consider the Eilenberg--Moore spectral sequence of the
bundle $p:B_TP\to BT^m$ with the fibre $\ZP$. This spectral sequence has
$E_2=\Tor_{k[v_1,\ldots,v_m]}\bigl(k(P),k\bigr)$, and $E_r\Rightarrow
H^{*}({\mathcal Z})$. It follows from inequality~(\ref{uneq}) that
$E_2=E_{\infty}$, which concludes the proof of the theorem.
\end{proof}

Let us turn again to the Eilenberg--Moore filtration
$\{F^{-p}H^{*}(\ZP)\}$ in $H^*(\ZP)$. The Poncar\'e duality defines a
filtration $\{F^{-p}H_{*}(\ZP)\}$ in the homology of $\ZP$.
It turns out that elements
from $F^{-1}H_{*}({\ZP})$ have very transparent geometric realization.
Namely, the following statement holds:

\begin{theorem}
\label{realiz}
  Elements of
  $\{F^{-1}H_{*}(\ZP)\}$ can be realized as embedded submanifolds of $\ZP$.
  These submanifolds are spheres of odd dimensions for
  generators of $\{F^{-1}H_{*}(\ZP,\Z)\}$
\end{theorem}
\begin{proof}
It follows from Theorem \ref{cohomZ} that
$$
  F^{-1}H^{*}(\ZP)/F^{0}H^{*}(\ZP)=
  \Tor^{-1}_{k[v_1,\ldots,v_m]}\bigl(k(P),k\bigr).
$$
By Proposition \ref{f0}, $F^{0}H^{*}(\ZP)=H^{0}(\ZP)$.
Take a basis of $\Tor^{-1}_{k[v_1,\ldots,v_m]}\bigl(k(P),k\bigr)$ consisting
of elements $v_{i_1\ldots i_p}$ of internal degree $2p$ such that
$v_{i_1},\ldots,v_{i_p}$ are primitive collections of the vertices of $K_P$
(see the proof of Theorem~\ref{cohomZ}).  If $v_{i_1},\ldots,v_{i_p}$ is a
primitive collection, then the subcomplex of $K_P$ consisting of all
simplices with vertices among $v_{i_1},\ldots,v_{i_p}$ is a simplicial
complex consisting of all faces of a simplex except one of the highest
dimension (i.e. the boundary of a simplex). In terms of the simple polytope
$P$ the element $v_{i_1\ldots i_p}$ corresponds to the set
$\{F_{i_1},\ldots,F_{i_p}\}$ of codimension-one faces such that
$F_{i_1}\cap\cdots\cap F_{i_p}=\emptyset$, though any proper subset of
$\{F_{i_1},\ldots,F_{i_p}\}$ has a non-empty intersection. Note that the
element $v_{i_1\ldots i_p}\in\Tor^{-1}_{k[v_1,\ldots,v_m]}\bigl(k(P),k\bigr)$
defines, by means of the isomorphism from Theorem~\ref{cohomZ}, an element of
$H^{*}(\ZP)$ of dimension $2p-1$.  Now, we take one point inside each face
$F_{i_1}\cap\cdots\cap\widehat{F_{i_r}}\cap\cdots\cap F_{i_p}$,
$1\le r\le p$, ($F_{i_r}$ is dropped); then we can embed the simplex
$\Delta^{p-1}$ on these points into the polytope $P$ in such a way that the
boundary $\partial\Delta^{p-1}$ embeds into $\partial P$. (Compare this with
the construction of the cubical decomposition of $P$ in Theorem~\ref{Pcub}.)
Let $\rho:\ZP=(T^m\times P^n)/\!\sim\:\to P^n$ be the projection
onto the orbit space; then it can be easily seen that
$\rho^{-1}(\Delta^{p-1})=(T^p\times\Delta^{p-1})/\!\sim)\times T^{m-p}=
S^{2p-1}\times T^{m-p}$. In this way
we obtain an embedding $S^{2p-1}\hookrightarrow\ZP$ that realize the
element of $H_{*}(\ZP)$ dual to $v_{i_1\ldots i_p}$.
\end{proof}

\subsection{Multiplicative structure of the cohomology of $\ZP$.}

Here we describe the ring $H^{*}(\ZP)$.

The bigraded $k$-module
$\Tor_{k[v_1,\ldots,v_m]}\bigr(k(P),k\bigl)$ can be calculated either by
means of a resolution of the face ring $k(P)$ or by means of a resolution
of~$k$. In the previous subsection we studied the minimal resolution of
$k(P)$ as a $k[v_1,\ldots,v_m]$-module. Here we use another approach based on
the {\it Koszul resolution} of $k$ as a $k[v_1,\ldots,v_m]$-module.
This allows to invest the bigraded $k$-module
$\Tor_{k[v_1,\ldots,v_m]}\bigl(k(P),k\bigr)$ with a bigraded $k$-algebra
structure. We show that the corresponding total graded $k$-algebra is
isomorphic to the algebra $H^{*}(\ZP)$. This approach also gives
us the description of $H^{*}(\ZP)$ as a cohomology algebra of some
differential (bi)graded algebra.

\medskip

Let $\Gamma=k[y_1,\ldots,y_n]$, $\deg y_i=2$, be a graded polynomial algebra
over $k$, and let $\Lambda[u_1,\ldots,u_n]$ denote an exterior algebra over
$k$ on generators $u_1,\ldots,u_n$. Consider the bigraded differential
algebra
$$
  \mathcal E=\Gamma\otimes\Lambda[u_1,\ldots,u_n],
$$
whose gradings and differential are defined by
$$
\begin{array}{lclclcl}
  \bideg(y_i\otimes1)&=&(0,2),&\quad&d(y_i\otimes1)&=&0;\\
  \bideg(1\otimes u_i)&=&(-1,2),&\quad&d(1\otimes u_i)&=&y_i\otimes1,
\end{array}
$$
and requiring that $d$ be a derivation of algebras.
The differential adds $(1,0)$ to bidegree, hence, the components
$\mathcal E^{-i,*}$ form a cochain complex. This complex will be also
denoted by $\mathcal E$. It is well known that this complex defines a
$\G$-free resolution of $k$ (regarded as a $\G$-module) called the
{\it Koszul resolution} (see~\cite{Ma}).

\begin{proposition}
\label{exter}
  Let $\Gamma=k[y_1,\ldots,y_n]$, and let $A$ be a $\Gamma$-module,
  then
  $$
    \Tor_{\Gamma}(A,k)=H\bigl[A\otimes\Lambda[u_1,\ldots,u_n],d\bigr],
  $$
  where $d$ is defined as $d(a\otimes u_i)=(y_i\cdot a)\otimes1$ for any
  $a\in A$.
\end{proposition}
\begin{proof}
Let us consider the introduced above $\Gamma$-free Koszul resolution
$\mathcal E=\G\otimes\Lambda[u_1,\ldots,u_n]$ of $k$. Then
$$
  \Tor_{\Gamma}(A,k)=H\bigl[A\otimes_{\Gamma}\Gamma\otimes
  \Lambda[u_1,\ldots,u_n],d\bigr]=
  H\bigl[A\otimes\Lambda[u_1,\ldots,u_n],d\bigr].
$$
\end{proof}

Now let us consider the principal $T^m$-bundle $\ZP\times ET^m\to B_TP$
pulled back from the universal $T^m$-bundle by the map $p:B_TP\to BT^m$
(see~(\ref{kill})). The following lemma holds.

\begin{figure}
\begin{center}
\begin{picture}(32,30)
  \multiput(0,0)(8,0){4}{\line(0,1){24}}
  \multiput(0,0)(0,8){3}{\line(1,0){32}}
  \put(3,3){$1$}
  \put(2,11){$u_i$}
  \put(18,3){$v_i$}
  \put(6,10){\vector(3,-1){13}}
  \put(10,9){\footnotesize $d_2$}
\end{picture}%
\caption{The $E_2$ term of the spectral sequence for
$\ZP\times ET^m\to B_TP$.}
\label{fig3}
\end{center}
\end{figure}
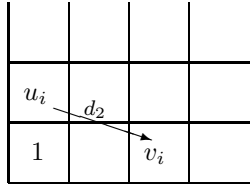

\begin{lemma}
\label{3termgen}
  The following isomorphism describes the $E^{(s)}_3$ term
  of the Leray--Serre spectral sequence $\{E^{(s)}_r,d_r\}$ of
  the bundle $\ZP\times ET^m\to B_TP$:
  $$
    E^{(s)}_3\cong\Tor_{k[v_1,\ldots,v_m]}\bigl(k(P),k\bigr).
  $$
\end{lemma}
\begin{proof}
First, consider the $E^{(s)}_2$ term of the
spectral sequence. Since
$H^*(T^m)=\Lambda[u_1,\ldots,u_m]$, $H^*(B_TP)=k(P)=k[v_1,\ldots,v_m]/I$,
we have
$$
  E^{(s)}_2=k(P)\otimes\Lambda[u_1,\ldots,u_m].
$$
It can be easily seen that the differential $d_2^{(s)}$ acts as follows:
$$
  d_2^{(s)}(1\otimes u_i)=v_i\otimes 1,\quad d_2^{(s)}(v_i\otimes 1)=0
$$
(see Figure~\ref{fig3}).
Now, since $E^{(s)}_3=H[E_2^{(s)},d_2^{(s)}]$, our assertion follows from
Proposition~\ref{exter} by putting $\G=k[v_1,\ldots,v_m]$, $A=k(P)$.
\end{proof}

Now we are ready to prove our main result on the cohomology of $\ZP$.

\begin{theorem}
\label{mult}
  The following isomorphism of graded algebras holds:
  \begin{gather*}
    H^{*}(\ZP)\cong H\bigl[k(P)\otimes\Lambda[u_1,\ldots,u_m],d\bigr],\\
    \bideg v_i=(0,2),\quad\bideg u_i=(-1,2),\\
    d(1\otimes u_i)=v_i\otimes 1,\quad d(v_i\otimes 1)=0.
  \end{gather*}
  Hence, the Leray--Serre spectral sequence of the $T^m$-bundle
  $\ZP\times ET^m\to B_TP$ collapses in the $E_3$ term.
\end{theorem}
\begin{proof}
Let us consider the bundle $p:B_TP\to BT^m$ with the fibre $\ZP$.
It follows from Theorem~\ref{cell} that the correspondent cochain algebras
are $C^{*}(BT^m)=k[v_1,\ldots,v_m]$ and $C^{*}(B_TP)=k(P)$, and the action of
$C^{*}(BT^m)$ on $C^{*}(B_TP)$ is defined by the quotient projection.
It was shown in~\cite[Proposition~3.4]{Sm} that there is an isomorphism
of algebras
$$
  \theta^{*}:\Tor_{C^{*}(BT^m)}\bigl(C^{*}(B_TP),k\bigr)\to H^{*}(\ZP).
$$
Now, it follows from above arguments and Proposition~\ref{exter} that
$$
  \Tor_{C^{*}(BT^m)}\bigl(C^{*}(B_TP),k\bigr)\cong
  H\bigl[k(P)\otimes\Lambda[u_1,\ldots,u_m],d\bigr],
$$
which concludes the proof.
\end{proof}

\subsection{Cohomology of $\ZP$ and torus actions.}
\label{addprop}

First, we consider the case where the simple polytope $P^n$ can be
realized as the orbit space for some quasitoric manifold (see
subsection~1.3). We show that the existence of a quasitoric manifold
enables us to reduce calculating of the cohomology of $\ZP$ to calculating
the cohomology of an algebra, which is much smaller than that from
Theorem~\ref{mult}.

As it was already discussed in section~1.3, a quasitoric manifold $M^{2n}$
over $P^n$ defines a principal $T^{m-n}$-bundle $\ZP\to M^{2n}$. This
bundle is induced from the universal $T^{m-n}$-bundle by a certain map
$f:M^{2n}\to BT^{m-n}$.

\begin{theorem}
\label{compar}
  Suppose $M^{2n}$ is a quasitoric manifold over a simple polytope
  $P^n$; then the Eilenberg--Moore spectral sequences of the
  commutative squares
  $$
  \begin{array}{ccccccc}
    \ZP\times ET^m&\longrightarrow&ET^m&&\ZP&\longrightarrow&ET^{m-n}\\
    \downarrow&&\downarrow&
    \qquad\text{and}\qquad&\downarrow&&\downarrow\\
    B_TP&\stackrel{p}{\longrightarrow}&BT^m&&
    M^{2n}&\stackrel{f}{\longrightarrow}&BT^{m-n}
  \end{array}
  $$
  are isomorphic.
\end{theorem}
\begin{proof}
Let $\{E_r,d_r\}$ be the Eilenberg--Moore spectral sequence of the first
commutative square, and let $\{\bar{E}_r,\bar{d}_r\}$ be that of the second
one.  Then, as it follows from the results of~\cite{EM},~\cite{Sm}, the
inclusions $BT^{m-n}\to BT^m$, $ET^{m-n}\to ET^m$, $M^{2n}\to B_TP$, and
$\ZP\to \ZP\times ET^m$ define a homomorphism of spectral sequences:
$g:\{E_r,d_r\}\to\{\bar{E}_r,\bar d_r\}$. First, we prove that
$g_2:E_2\to \bar E_2$ is an isomorphism.

The map $H^{*}(BT^m)\to H^*(BT^{m-n})$ is the quotient
projection $k[v_1,\ldots,v_m]\to k[w_1,\ldots,w_{m-n}]$ with the kernel
$J=(\l_1,\ldots,\l_n)$. By Theorem~\ref{cohomM}, we have
$H^*(M^{2n})=k[v_1,\ldots,v_m]\,/\,I{+}J$. Hence,
$f^*:H^*(BT^{m-n})=k[v_1,\ldots,v_m]/J\to
k[v_1,\ldots,v_m]\,/\,I{+}J=H^*(M^{2n})$ is the quotient epimorphism.

The $E_2$-terms of our spectral sequences are
$E_2=\Tor_{k[v_1,\ldots,v_m]}\bigl(k(P),k\bigr)$ and
$\bar{E}_2=\Tor_{k[w_1,\ldots,w_{m-n}]}\bigl(k(P)/J,k\bigr)$.

To proceed further we need the following result.

\begin{proposition}
\label{change}
  Let $\Lambda$ be an algebra and $\Gamma$ a subalgebra and set
  $\Omega=\Lambda//\Gamma$. Suppose that $\Lambda$ is a free $\Gamma$-module
  and we are given a right $\Omega$-module $A$ and a left $\Lambda$-module
  $C$. Then there exists a spectral sequence $\{E_r,d_r\}$ with
  $$
    E_r\Rightarrow \Tor_{\Lambda}(A,C),\quad
    E_2^{p,q}=\Tor^p_{\Omega}\bigl(A,\Tor_{\Gamma}^q(C,k)\bigr).
  $$
\end{proposition}
\begin{proof}
See \cite[p.349]{CE}.
\end{proof}

The next proposition is a modification of one assertion from~\cite{Sm}.

\begin{proposition}
\label{tortor}
  Suppose $f:k[v_1,\ldots,v_m]\to A$ is an epimorphism of graded algebras,
  $\deg v_i=2$, and $J\subset A$ is an ideal
  generated by a length $n$ regular sequence of degree-two elements of $A$.
  Then the following isomorphism holds:
  $$
    \Tor_{k[v_1,\ldots,v_m]}(A,k)=\Tor_{k[w_1,\ldots,w_{m-n}]}(A/J,k).
  $$
\end{proposition}
\begin{proof}
Let $J=(\l_1,\ldots,\l_n)$, $\deg\l_i=2$, and
$\{\l_1,\ldots,\l_n\}$ is a regular sequence.
Let $\hat\l_i$, $1\le i\le n$, be degree-two elements of
$k[v_1,\ldots,v_m]$ such that
$f(\hat\l_i)=\l_i$. Hence, $\hat\l_i=\l_{i1}v_1+\ldots+\l_{im}v_m$ and
$\mathop{\rm rk}(\l_{ij})=n$. Let us take elements $w_1,\ldots,w_{m-n}$
of degree two such that
$$
  k[v_1,\ldots,v_m]=k[\hat\l_1,\ldots,\hat\l_n,w_1,\ldots,w_{m-n}],
$$
and put $\Gamma=k[\hat\l_1,\ldots,\hat\l_n]$. Then $k[v_1,\ldots,v_m]$ is
a free $\Gamma$-module, and therefore,
by Proposition~\ref{change}, we have a spectral sequence
$$
  E_r\Rightarrow \Tor_{k[v_1,\ldots,v_m]}(A,k),\quad
  E_2=\Tor_{\Omega}\bigl(\Tor_{\Gamma}(A,k),k\bigr),
$$
where $\Omega=k[v_1,\ldots,v_m]//{\Gamma}=k[w_1,\ldots,w_{m-n}]$.

Since $\l_1,\ldots,\l_n$ is a regular sequence, $A$ is a free $\Gamma$-module.
Therefore,
\begin{align*}
  \Tor_{\Gamma}(A,k)&=A\otimes_{\Gamma}k=A/J\text{\quad and \quad}
  \Tor_{\Gamma}^q(A,k)=0\text{\quad for }\;q\ne 0,\\
  &\Rightarrow \quad E_2^{p,q}=0\text{\quad for }\;q\ne 0,&{}\\
  &\Rightarrow \quad \Tor_{k[v_1,\ldots,v_m]}(A,k)=
  \Tor_{k[w_1,\ldots,w_{m-n}]}(A/J,k),&{}
\end{align*}
which concludes the proof of the proposition.
\end{proof}

Now, we return to the proof of Theorem~\ref{compar}.
Setting $A=k(P)$ in Proposition~\ref{tortor} we deduce that
$g_2:E_2\to \bar E_2$ is an isomorphism.
The $E_2$ terms of both spectral sequences contain
only finite number of non-zero modules. In this situation
a homomorphism $g$ that defines an isomorphism in the
$E_2$ terms is an isomorphism of the spectral sequences
(see~\cite[XI, Theorem~1.1]{Ma}). Thus, Theorem~\ref{compar} is proved.
\end{proof}

\begin{corollary}
\label{cZ2}
  Suppose that $M^{2n}$ is a quasitoric manifold over a
  simple polytope $P^n$. Then the cohomology of $\ZP$ can be calculated as
  $$
    H^{*}(\ZP)=\Tor_{k[w_1,\ldots,w_{m-n}]}\bigl(H^{*}(M^{2n}),k\bigr).
  $$
\end{corollary}
\begin{proof}
By Theorem~\ref{mult},
$H^{*}(\ZP)=\Tor_{k[v_1,\ldots,v_m]}\bigl(k(P),k\bigr)$. Hence, our
assertion follows from the isomorphism between the $E_2$ terms of the
spectral sequences from Theorem~\ref{compar}.
\end{proof}

Let us turn again to the principal $T^{m-n}$-bundle $\ZP\to M^{2n}$ defined
by a quasitoric manifold $M^{2n}$. The following lemma is
analogous to Lemma~\ref{3termgen} (and is proved similarly).
\begin{lemma}
\label{3term}
  The following isomorphism holds for
  the Leray--Serre spectral sequence of the bundle $\ZP\to M^{2n}$:
  $$
    E^{(s)}_3\cong\Tor_{k[w_1,\ldots,w_{m-n}]}\bigl(H^*(M^{2n}),k\bigr)=
    \Tor_{k[w_1,\ldots,w_{m-n}]}(k(P)/\!J\:,k),
  $$
  where $E^{(s)}_3$ is the $E_3$ term of the Leray--Serre spectral sequence,
  and $H^*(M^{2n})\cong k(P)/\!J$ is invested with
  a $k[w_1,\ldots,w_{m-n}]$-module structure by means of the map
  $$
    k[w_1,\ldots,w_{m-n}]=k[v_1,\ldots,v_m]/J\to
    k[v_1,\ldots,v_m]/I{+}J= H^*(M^{2n}).\quad\square
  $$
\end{lemma}

\begin{theorem}
\label{degene3}
  Suppose $M^{2n}$ is a quasitoric manifold over $P^n$. Then the
  Leray--Serre spectral sequence of the principle $T^{m-n}$-bundle
  $\ZP\to M^{2n}$ collapses in the $E_3$ term, i.e. $E_3=E_{\infty}$.
  Furthermore, the following isomorphism of algebras holds
  \begin{gather*}
    H^{*}(\ZP)=H\bigl[(k(P)/\!J)\otimes
    \Lambda[u_1,\ldots,u_{m-n}],d\bigr],\\
    \bideg a=(0,\deg a),\quad\bideg u_i=(-1,2);\\
    d(1\otimes u_i)=w_i\otimes 1,\quad d(a\otimes 1)=0,
  \end{gather*}
  where $a\in k(P)/\!J=k[w_1,\ldots,w_{m-n}]/I$ and
  $\Lambda[u_1,\ldots,u_{m-n}]$ is an exterior algebra.
\end{theorem}
\begin{proof}
The cohomology algebra
$H\bigl[(k(P)/\!J)\otimes\Lambda[u_1,\ldots,u_{m-n}],d\bigr]$ is exactly
the $E_3$ term of the Leray--Serre spectral sequence for the bundle
$\ZP\to M^{2n}$. At the same time, it follows from Proposition~\ref{exter}
that this cohomology algebra is isomorphic to
$\Tor_{k[w_1,\ldots,w_{m-n}]}\bigl(H^{*}(M^{2n}),k\bigr)$.
Corollary~\ref{cZ2} shows that this is exactly $H^{*}(\ZP)$. Since the
Leray--Serre spectral sequence converges to $H^{*}(\ZP)$, it follows that
it collapses in the $E_3$ term.
\end{proof}

The algebra $\bigl(k(P)/\!J\bigr)\otimes\Lambda[u_1,\ldots,u_{m-n}]$ from
Theorem~\ref{degene3} is significantly smaller than the algebra
$k(P)\otimes\Lambda[u_1,\ldots,u_m]$ from general Theorem~\ref{mult}.
This enables to calculate the cohomology of $\ZP$ more efficiently.

\medskip

A rank $m-n$ torus subgroup of $T^m$ that acts freely on $\ZP$ gives rise to
a quasitoric manifold $M^{2n}=\ZP/T^{m-n}$ with orbit space $P^n$.
In the general case, such a subgroup may fail to exist; however, one still
may be able to find a subgroup of dimension less than $m-n$ that acts
freely on $\ZP$. So, suppose that a subgroup $H\cong T^r$ acts on $\ZP$
freely. Then the inclusion $s:H\hookrightarrow T^m$ is defined by an
integer $(m\times r)$-matrix $S=(s_{ij})$ such that the $\Z$-module
spanned by its columns $s_j=(s_{1j},\ldots,s_{mj})^\top$, $j=1,\ldots,r$ is
a direct summand in $\Z^m$. Choose any basis $t_i=(t_{i1},\ldots,t_{im})$,
$i=1,\ldots,m-r$ in the kernel of the dual map
$s^{*}:(\Z^m)^{*}\to(\Z^r)^{*}$. Then the cohomology ring of the quotient
manifold ${\mathcal Y}_{(r)}=\ZP/H$ is described by the following theorem,
which generalize both Corollary~\ref{Mtor} and
Theorem~\ref{cohomZ}.

\begin{theorem}
\label{quot}
  The following isomorphism of algebras holds:
  $$
    H^{*}({\mathcal Y}_{(r)})\cong\Tor_{k[t_1,\ldots,t_{m-r}]}
    \bigl(k(P),k\bigr),
  $$
  where the action of $k[t_1,\ldots,t_{m-r}]$ on
  $k(P)=k[v_1,\ldots,v_m]/I$ is defined by the map
  $$
  \begin{array}{rcl}
    k[t_1,\ldots,t_{m-r}]&\to&k[v_1,\ldots,v_m]\\[1mm]
    t_i&\to&t_{i1}v_1+\ldots+t_{im}v_m.
  \end{array}
  $$
\end{theorem}

\begin{remark}
  Corollary~\ref{Mtor} corresponds to the value $r=m-n$, while
  Theorem~\ref{cohomZ} corresponds to the value $r=0$.
\end{remark}

\begin{proof}
The inclusion of the subgroup $H\cong T^r\to T^m$ defines a
map of classifying spaces $h:BT^r\to BT^m$. Let us consider the bundle
pulled back by this map from the bundle $p:B_TP\to BT^m$ with the fibre $\ZP$.
It follows directly from the construction of $B_TP$ (see subsection~1.2) that
the total space of this bundle has homotopy type
$\mathcal Y_{(r)}$ (more precisely, it is homeomorphic to
${\mathcal Y}_{(r)}\times ET^r$). Hence, we have the commutative square
$$
  \begin{array}{ccc}
    {\mathcal Y}_{(r)}&\longrightarrow&B_TP\\
    \downarrow&&\downarrow\\
    {BT^r}&\longrightarrow & BT^m.
  \end{array}
$$
The corresponding Eilenberg--Moore spectral sequence converges to the
cohomology of $\mathcal Y_{(r)}$ and has the following $E_2$ term:
$$
  E_2=\Tor_{k[v_1,\ldots,v_m]}\bigl(k(P),k[w_1,\ldots,w_r]\bigr),
$$
where the action of $k[v_1,\ldots,v_m]$ on $k[w_1,\ldots,w_r]$ is defined by
the map $s^{*}$, i.e. $v_i\to s_{i1}w_1+\ldots+s_{ir}w_r$.
Using~\cite[Proposition~3.4]{Sm} in the similar way as in the proof of
Theorem~\ref{mult}, we show that the spectral sequence collapses in the
$E_2$ term and the following isomorphism of algebras holds:
\begin{equation}
\label{Ytor1}
  H^{*}({\mathcal Y}_{(r)})=\Tor_{k[v_1,\ldots,v_m]}
  \bigl(k(P),k[w_1,\ldots,w_r]\bigr).
\end{equation}
Now put $\Lambda=k[v_1,\ldots,v_m]$,
$\Gamma=k[t_1,\ldots,t_{m-r}]$, $A=k[w_1,\ldots,w_r]$, and $C=k(P)$ in
Proposition~\ref{change}. Since $\Lambda$ here is a free
$\Gamma$-module and $\Omega=\Lambda//\Gamma=k[w_1,\ldots,w_r]$,
a spectral sequence $\{E_s,d_s\}$ arises. Its $E_2$ term is
$$
  E_2^{p,q}=\Tor^p_{k[w_1,\ldots,w_r]}
  \Bigl(A,\Tor^q_{k[t_1,\ldots,t_{m-r}]}\bigl(k(P),k\bigr)\Bigr),
$$
and it converges to
$\Tor_{k[v_1,\ldots,v_m]}\bigl(k(P),k[w_1,\ldots,w_r]\bigr)$.
Since $A$ is a free module over $k[w_1,\ldots,w_r]$
with one generator 1, we have
$$
  E_2^{p,q}=0\;\mbox{ for }p\ne0,\quad
  E_2^{0,q}=\Tor^q_{k[t_1,\ldots,t_{m-r}]}\bigl(k(P),k\bigr).
$$
Thus, the spectral sequence collapses in the $E_2$ term, and we have
the isomorphism of algebras:
$$
  \Tor_{k[v_1,\ldots,v_m]}\bigl(k(P),k[w_1,\ldots,w_r]\bigr)\cong
  \Tor_{k[t_1,\ldots,t_{m-r}]}\bigl(k(P),k\bigr),
$$
which together with the isomorphism~(\ref{Ytor1}) proves the theorem.
\end{proof}

Below we characterize subgroups $H\subset T^m$ that act on $\ZP$ freely.

Let us consider again the integer $(m\times r)$-matrix $S$ defining the
subgroup $H\subset T^m$ of rank $r$. For each vertex
$v=F_{i_1}\cap\cdots\cap F_{i_n}$ of the polytope $P^n$ denote by
$S_{i_1,\ldots,i_n}$ the $(m-n)\times r$-submatrix of $S$ that is obtained
by deleting the rows $i_1,\ldots,i_n$. In this way we construct $f_{n-1}$
submatrices of the size $(m-n)\times r$. Then the following criterion
for the freeness of the action of $H$ on $\ZP$ holds.

\begin{lemma}
\label{free}
  The action of the subgroup $H\subset T^m$ defined by an integer
  $(m\times r)$-matrix $S$ on the manifold $\ZP$ is free if and only
  if for any vertex $v=F_{i_1}\cap\ldots\cap F_{i_n}$ of $P^n$ the
  corresponding $(m-n)\times r$-submatrix $S_{i_1,\ldots,i_n}$ defines
  a direct summand $\Z^r\subset\Z^{m-n}$.
\end{lemma}
\begin{proof}
It follows from Definition \ref{defzp} that the orbits of the
action of $T^m$ on $\ZP$ corresponding to the vertices
$v=F_{i_1}\cap\ldots\cap F_{i_n}$ of $P^n$ have maximal (rank $n$)
isotropy subgroups. These isotropy subgroups are the coordinate
subgroups $T^n_{i_1,\ldots,i_n}\subset T^m$. A subgroup $H$ acts freely on
$\ZP$ if and only if it has only unit in the intersection with each isotropy
subgroup. This means that the $m\times(r+n)$-matrix obtained by adding
$n$ columns $(0,\ldots,0,1,0,\ldots,0)^\top$ (1 stands on the place $i_j$,
$j=1,\ldots,n$) to $S$ defines a direct summand $\Z^{k+n}\subset\Z^m$.
(This matrix corresponds to the subgroup
$H\times T^n_{i_1,\ldots,i_n}\subset T^m$.) Obviously, this
is equivalent to the requirements of the lemma.
\end{proof}

In particular, for subgroups of rank $m-n$ we obtain

\begin{corollary}
\label{maxfree}
  The action of the rank $m-n$ subgroup $H\subset T^m$ defined by an integer
  $m\times(m-n)$-matrix $S$ on the manifold $\ZP$ is free if and only
  if for any vertex $v=F_{i_1}\cap\ldots\cap F_{i_n}$ of $P^n$
  the minor $(m-n)\times(m-n)$-matrix $S_{i_1\ldots i_n}$
  has $\det S_{i_1\ldots i_n}=\pm1$.\ep
\end{corollary}

\begin{remark}
  Compare this with Proposition~\ref{fprop} and Theorem~\ref{zu}.
  Note that unlike the situation of Theorem~\ref{zu}, the subgroup
  $H\cong T^{m-n}$ satisfying the condition of Corollary~\ref{maxfree} may
  fail to exist.
\end{remark}

The inclusion $s:\Z^{m-n}\to\Z^m$ defines the short exact sequence
$$
\begin{CD}
  0 @>>> \Z^{m-n} @>s>> \Z^m @>>> \Z^n @>>> 0.
\end{CD}
$$
It can be easily seen that the condition from Corollary~\ref{maxfree}
is equivalent to the following: the map $\Z^m\to\Z^n$ above
is a characteristic
function in the sense of Definition~\ref{chf}. Thus, we have obtained the
another interpretation of the fact that quasitoric manifolds exist over
$P^n$ if and only if it is possible to find a subgroup $H\cong T^{m-n}$
that acts on $\ZP$ freely.

\medskip

As it follows from Lemma \ref{free}, the one-dimensional subgroup
$H\cong T^1$ defined to the diagonal inclusion $T^1\subset T^m$
always acts on $\ZP$ freely. Indeed, in this situation the matrix $S$ is a
column of $m$ units and the condition from Lemma~\ref{free} is obviously
satisfied.
Theorem~\ref{quot} gives the following formula for the cohomology of
the corresponding quotient manifold ${\mathcal Y}_{(1)}=\ZP/H$:
\begin{equation}
\label{y1}
  H^{*}({\mathcal Y}_{(1)})\cong\Tor_{k[t_1,\ldots,t_{m-1}]}
  \bigl(k(P),k\bigr),
\end{equation}
where the action of $k[t_1,\ldots,t_{m-1}]$ on $k(P)=k[v_1,\ldots,v_m]/I$ is
defined by the homomorphism
$$
  \begin{array}{rcl}
    k[t_1,\ldots,t_{m-1}]&\to&k[v_1,\ldots,v_m],\\[1mm]
    t_i&\to&v_i-v_m.
  \end{array}
$$
The principal $T^1$-bundle $\ZP\to{\mathcal Y}_{(1)}$ is pulled back from the
universal $T^1$-bundle by a certain map
$c:{\mathcal Y}_{(1)}\to BT^1=\C P^{\infty}$. Since $H^{*}(\C
P^{\infty})=k[v]$, $v\in H^2(\C P^{\infty})$, the element $c^{*}(v)\in
H^2({\mathcal Y}_{(1)})$ is defined. Then, the following statement holds.

\begin{lemma}
\label{neib}
  A polytope $P^n$ is $q$-neighbourly if and only if
  $\bigl(c^{*}(v)\bigr)^q\ne0$.
\end{lemma}
\begin{proof}
The map $c^{*}$ takes the cohomology ring $k[v]$ of
$\C P^{\infty}$ to the subring
$k(P)\otimes_{k[t_1,\ldots,t_{m-1}]}k=
\Tor^0_{k[t_1,\ldots,t_{m-1}]}\bigl(k(P),k\bigr)$
of the cohomology ring of ${\mathcal Y}_{(1)}$ (see~(\ref{y1})). This subring is
isomorphic to the quotient ring $k(P)/(v_1=\ldots=v_m)$. Now, the
assertion follows from the fact that a polytope $P^n$ is
$q$-neighbourly if and only if the ideal $I$ (see Definition~\ref{frpol})
does not contain monomials of degree less than $q+1$.
\end{proof}

Now we return to the general case of a subgroup $H\cong T^r$ acting
on $\ZP$ freely. For such a subgroup we have
\begin{multline*}
  B_TP=\ZP\times_{T^m}ET^m=\left((\ZP/T^r)\times_{T^{m-r}}ET^{m-r}\right)
  \times ET^r\\
  =({\mathcal Y}_{(r)}\times_{T^{m-r}}ET^{m-r})\times ET^r.
\end{multline*}
Hence, there is defined a principal $T^{m-r}$-bundle
${\mathcal Y}_{(r)}\times ET^m\to B_TP$.

\begin{theorem}
\label{degene3gen}
  The Leray--Serre spectral sequence of the $T^{m-r}$-bundle
  ${\mathcal Y}_{(r)}\times ET^m\to B_TP$ collapses in the
  $E_3$ term, i.e. $E_3=E_{\infty}$. Furthermore,
  \begin{gather*}
    H^{*}({\mathcal Y}_{(r)})=H\bigl[k(P)\otimes
    \Lambda[u_1,\ldots,u_{m-r}],d\bigr],\\
    d(1\otimes u_i)=(t_{i1}v_1+\ldots+t_{im}v_m)\otimes 1,
    \quad d(a\otimes 1)=0;\\
    \bideg a=(0,\deg a),\quad\bideg u_i=(-1,2),
  \end{gather*}
  where $a\in k(P)=k[v_1,\ldots,v_m]/I$ and $\Lambda[u_1,\ldots,u_{m-r}]$ is
  an exterior algebra.
\end{theorem}
\begin{proof}
In the similar way as in Lemma~\ref{3termgen} we show that
the $E_3$ term of the spectral sequence is
$$
  E_3=H\bigl[k(P)\otimes\Lambda[u_1,\ldots,u_{m-k}],d\bigr]=
  \Tor_{k[t_1,\ldots,t_{m-r}]}\bigl(k(P),k\bigr).
$$
Theorem~\ref{quot} shows that this is exactly
$H^{*}({\mathcal Y}_{(r)})$.
\end{proof}

\begin{remark}
  Theorem \ref{mult} and Corollary \ref{Mtor} can be obtained from
  this theorem by setting $r=0$ and $r=m-n$ respectively.
\end{remark}

\subsection{Explicit calculation of $H^*(\ZP)$ for some particular
polytopes.}

\ \\
1. Our first example demonstrates how the above methods work in the simple
case when $P$ is a product of simplices. So, let
$P^n=\D^{i_1}\times\D^{i_2}\times\ldots\times\D^{i_k}$, where $\D^i$ is an
$i$-simplex and $\sum_ki_k=n$. This $P^n$ has $n+k$ facets, i.e.
$m=n+k$. Lemma~\ref{prod} shows that
$\ZP={\mathcal Z}_{\D^{i_1}}\times\ldots\times{\mathcal Z}_{\D^{i_k}}$.

The minimal resolution (\ref{resol}) of $k(P_i)$ in the case $P_i=\D^i$ is as
follows
$$
  0\longrightarrow R^{-1}\stackrel{d^{-1}}{\longrightarrow}
  R^{0}\stackrel{d^{0}}{\longrightarrow}k(P_i)\longrightarrow 0,
$$
where $R^0$, $R^{-1}$ are free one-dimensional
$k[v_1,\ldots,v_{i+1}]$-modules and $d^{-1}$ is the multiplication by
$v_{1}\cdot\ldots\cdot v_{i+1}$. Hence, we have the isomorphism of algebras
$$
  \Tor_{k[v_1,\ldots,v_{i+1}]}\bigl(k(P_i),k\bigr)=
  \Lambda[a],\quad\bideg a=(-1,2i+2),
$$
where $\Lambda[a]$ is an exterior $k$-algebra on one generator $a$. Now,
Theorem~\ref{mult} shows that
$$
  H^{*}({\mathcal Z}_{\D^i})=\Lambda[a],\quad\deg a=2i+1,
$$
Thus, the cohomology of
$\ZP={\mathcal Z}_{\D^{i_1}}\times\ldots\times{\mathcal Z}_{\D^{i_k}}$ is
$$
  H^{*}(\ZP)=\Lambda[a_1,\ldots,a_k],\quad\deg a_l=2i_l+1.
$$

Actually, Example \ref{sphere} shows that our $\ZP$ is the product of
spheres: $\ZP=S^{2i_1+1}\times\ldots\times S^{2i_k+1}$. However, our
calculation of the cohomology does not use the geometrical constructions from
section~2.

\medskip

2. In our next example we consider plane polygons,
i.e. the case $n=2$. Let $P^2$ be a convex $m$-gon. Then the corresponding
manifold $\ZP$ is of dimension $m+2$. First, we compute the Betti numbers of
these manifolds.

It can be easily seen that there is at least one quasitoric manifold $M^4$
over $P^2$.
Let us consider the $E_2$ term of the Leray--Serre spectral sequence for the
bundle $q:\ZP\to M^4$ with the fibre $T^{m-2}$.
Theorem~\ref{cohomM} shows that
$H^2(M^{4})$ has rank $m-2$ and the ring
$H^{*}(M^{4})$ is multiplicativelly generated by elements of degree~2.
The ring $H^{*}(T^{m-2})$ is an exterior algebra.
We choose bases $w_1,\ldots,w_{m-2}$ in $H^2(M^{4})$ and
$u_1,\ldots,u_{m-2}$ in $H^2(T^{m-2})$ such that
the second differential of the spectral sequence takes $u_i$ to $w_i$
(more precisely, $d_2(u_i\otimes 1)=1\otimes w_i$, see Figure~4).
Furthermore, the map
$q^{*}:H^{*}(M^{4})\to H^{*}(\ZP)$ is zero homomorphism in degrees $\ge0$.
This follows from the fact that the map $f^{*}:H^{*}(BT^{m-2})\to H^{*}(M^4)$
is epimorphic (see the proof of Theorem~\ref{compar}) and from
the commutative diagram
$$
  \begin{CD}
    \ZP @>>> ET^{m-2}\\
    @VqVV @VVV\\
    M^4 @>f>> BT^{m-2}.
  \end{CD}
$$
Using all these facts and Corollary~\ref{degene3} (which gives
$E_3=E_{\infty}$), we deduce
that all differentials $d_2^{0,*}$ are monomorphisms, and all
differentials $d_2^{2,*}$ are epimorphisms.

\begin{figure}
\begin{center}
\begin{picture}(32,60)
  \multiput(3,5)(8,0){6}{\line(0,1){26}}
  \multiput(3,5)(0,8){7}{\line(1,0){42}}
  \multiput(3,35)(8,0){6}{\line(0,1){20}}
  \put(6,8){$1$}
  \put(5,16){$u_i$}
  \put(21,8){$w_i$}
  \multiput(9,15)(16,0){2}{\vector(3,-1){13}}
  \multiput(9,23)(16,0){2}{\vector(3,-1){13}}
  \multiput(9,47)(16,0){2}{\vector(3,-1){13}}
  \put(12.5,22.5){\footnotesize $d_2^{0,*}$}
  \put(28.2,22.5){\footnotesize $d_2^{2,*}$}
  \put(6.5,1.5){\footnotesize 0}
  \put(22.5,1.5){\footnotesize 2}
  \put(38.5,1.5){\footnotesize 4}
  \put(0,8){\footnotesize 0}
  \put(2,48){\llap{\footnotesize $m{-}2$}}
  \put(4,33){\hbox to 38mm{\dotfill}}
\end{picture}%
\caption{The $E_2$ term of the spectral sequence for
$q:\ZP\to M^4$.}
\label{fig4}
\end{center}
\end{figure}
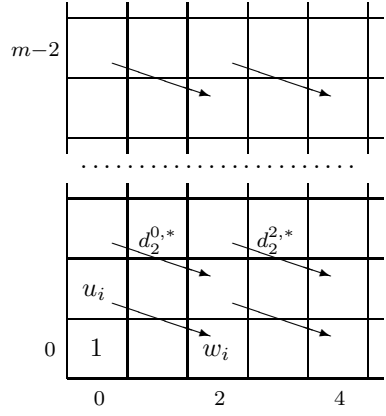

Now, using Theorem~\ref{degene3} we obtain by easy calculations the following
formulae for the Betti numbers $b^i(\ZP)$:
\begin{align}
  &b^0(\ZP)=b^{m+2}(\ZP)=1,\notag\\
  &b^{1}(\ZP)=b^{2}(\ZP)=b^{m}(\ZP)=
  b^{m+1}(\ZP)=0,\notag\\
  &b^{k}(\ZP)=(m-2)\binom{m-2}{k-2}-\binom{m-2}{k-1}-
  \binom{m-2}{k-3}\\
  &\hspace{0.3\textwidth}=\binom{m-2}{k-3}\frac{m(m-k)}{k-1},
  \quad3\le k\le m-1.
  \notag
\end{align}
For small $m$ the above formulae give us the following:
$$
\begin{array}{ll}
  m=3:& b^0({\mathcal Z}^5)=b^5({\mathcal Z}^5)=1;\\[1mm]
  m=4:& b^0({\mathcal Z}^6)=b^6({\mathcal Z}^6)=1,\quad
  b^3({\mathcal Z}^6)=2,
\end{array}
$$
(all other Betti numbers are zero).
Both cases are covered by the previous example. Indeed, for $m=3$ we have
$P^2=\D^2$, and for $m=4$ we have $P^2=\D^1\times\D^1$. As it was pointed out
above, in this cases $\ZP^5=S^5$, $\ZP^6=S^3\times S^3$. Further,
$$
\begin{array}{ll}
  m=5:& b^0({\mathcal Z}^7)=b^7({\mathcal Z}^7)=1,\quad
  b^3({\mathcal Z}^7)=b^4({\mathcal Z}^7)=5;\\[1mm]
  m=6:& b^0({\mathcal Z}^8)=b^8({\mathcal Z}^8)=1,\quad
  b^3({\mathcal Z}^8)=b^5({\mathcal Z}^8)=9,\quad
  b^4({\mathcal Z}^8)=16,
\end{array}
$$
(all other Betti numbers are zero), and so on.

Now we want to describe the ring structure in the cohomology.
Theorem~\ref{mult} gives us the isomorphism of algebras
\begin{equation}
\label{pg}
  H^{*}(\ZP^{m+2})\cong\Tor_{k[v_1,\ldots,v_m]}\bigl(k(P^2),k\bigr)
  =H\bigl[k(P^2)\otimes\Lambda[u_1,\ldots,u_m],d\bigr].
\end{equation}
If $m=3$, then $k(P)=k[v_1,v_2,v_3]/v_1v_2v_3$; if $m>3$ we have
$k(P)=k[v_1,\ldots,v_m]/I$, where $I$ is generated by monomials
$v_iv_j$ such that $i\ne j\pm1$. (Here we use the agreement
$v_{m+i}=v_i$ and $v_{i-m}=v_i$.) Below we give the complete description of
the multiplication in the case $m=5$. The general case is similar but more
involved. It is easy to check that five generators of
$H^3(\ZP)$ are represented by the cocycles
$v_i\otimes u_{i+2}\in
k(P^2)\otimes\Lambda[u_1,\ldots,u_m]$, $i=1,\ldots,5$, while five generators
of $H^4(\ZP)$ are represented by the cocycles $v_j\otimes u_{j+2}u_{j+3}$,
$j=1,\ldots,5$. The product of cocycles $v_i\otimes u_{i+2}$ and $v_j\otimes
u_{j+2}u_{j+3}$ represents a non-trivial cohomology class in $H^7(\ZP)$ if
and only if the set $\{i,i+2,j,j+2,j+3\}$ is the whole index set
$\{1,2,3,4,5\}$.  Hence, for each cohomology class $[v_i\otimes u_{i+2}]$
there is a unique (Poincar\'e dual) cohomology class $[v_j\otimes
u_{j+2}u_{j+3}]$ such that the product $[v_i\otimes u_{i+2}]\cdot[v_j\otimes
u_{j+2}u_{j+3}]$ is non-trivial. This product defines a fundamental
cohomology class of $\ZP$ (for example, it is represented by the cocycle
$v_1v_2\otimes u_3u_4u_5$). In the next section we prove the similar
statement in the general case. All other products in the cohomology algebra
$H^{*}(\ZP^7)$ are trivial.

\section{Cohomology of $\ZP$ and combinatorics of simple polytopes}

Theorem \ref{mult} shows that the cohomology of $\ZP$ is naturally a
bigraded algebra. The Poincar\'e duality in $H^{*}(\ZP)$ regards this
bigraded structure. More precisely, the Poincar\'e duality has the
following combinatorial interpretation.

\begin{lemma}
\label{pd}
In the bigraded differential algebra
$\bigl[k(P)\otimes\Lambda[u_1,\ldots,u_m],d\bigr]$ from
{\rm Theorem~\ref{mult}}
\begin{enumerate}
  \item For each vertex $v=F^{n-1}_{i_1}\cap\cdots\cap F^{n-1}_{i_n}$ of
  the polytope $P^n$
  the element $v_{i_1}\cdots v_{i_n}\otimes u_{j_1}\cdots u_{j_{m-n}}$,
  where $j_1<\ldots<j_{m-n}$,
  $\{i_1,\ldots,i_n,j_1,\ldots,j_{m-n}\}=\{1,\ldots,m\}$,
  represents the fundamental class of $\ZP$.
  \item Two cocycles
  $v_{i_1}\cdots v_{i_p}\otimes u_{j_1}\cdots u_{j_r}$ and
  $v_{k_1}\cdots v_{k_s}\otimes u_{l_1}\cdots u_{l_t}$ represent
  Poincar\'e dual cohomology classes in $H^{*}(\ZP)$ if and only if
  $p+s=n$, $r+t=m-n$, $\{i_1,\ldots,i_p,k_1,\ldots,k_s\}$ is the index set
  of facets meeting in some vertex $v\in P^n$, and
  $\{i_1,\ldots,i_p,j_1,\ldots,j_r,k_1,\ldots,k_s,l_1,\ldots,l_t\}
  =\{1,\ldots,m\}$.
\end{enumerate}
\end{lemma}
\begin{proof}
The first assertion follows from the fact that the cohomology
class under consideration is a generator of the module
$\Tor^{-(m-n),2m}_{k[v_1,\ldots,v_m]}
\bigl(k(P),k\bigr)\cong H^{m+n}(\ZP^{m+n})$
(see Theorem~\ref{mult}). The second
assertion holds since two cohomology classes are Poincar'e dual if
and only if their product is the fundamental cohomology class.
\end{proof}

In what follows we use the following notations:
$\T^{i}=\Tor^{-i}_{k[v_1,\ldots,v_m]}\bigl(k(P),k\bigr)$ and
$\T^{i,2j}=\Tor^{-i,2j}_{k[v_1,\ldots,v_m]}
\bigl(k(P),k\bigr)$. We define the
{\it bigraded Betti numbers} of $\ZP$ as
\begin{equation}
\label{bb}
  b^{-i,2j}(\ZP)=\dim_k\Tor^{-i,2j}_{k[v_1,\ldots,v_m]}
  \bigl(k(P),k\bigr).
\end{equation}
Then Theorem \ref{cohomZ} can be reformulated as
$b^k(\ZP)=\sum_{2j-i=k}b^{-i,2j}(\ZP)$.
The second part of Lemma~\ref{pd} shows that
$b^{-i,2j}(\ZP)=b^{-(m-n-i),2(m-j)}(\ZP)$ for all $i,j$. These equalities
can be written as the following identities for the Poincar\'e series
$F(\T^i,t)=\sum_{r=0}^mb^{-i,2r}t^{2r}$ of $\T^i$:
\begin{equation}
\label{algdual}
  F(\T^i,t)=t^{2m}F\left(\T^{m-n-i},\mbox{$\frac 1t$}\right),\quad
  i=1,\ldots,m-n.
\end{equation}
It is well known in commutative algebra that the above identities hold for
the so-called Gorenstein rings (see~\cite{St}). The face ring of a simplicial
subdivision of sphere is a Gorenstein ring. In particular, the ring
$k(P^n)$ is Gorenstein for any simple polytope $P^n$.

A simple combinatorial argument (see~\cite[part~II, \S1]{St})
shows that for any $(n-1)$-dimensional simplicial complex $K$ the
Poincar\'e series $F\bigl(k(K),t\bigr)$ can be written as follows
$$
 F\bigl(k(K),t\bigr)=1+\sum_{i=0}^{n-1}\frac{f_i t^{2(i+1)}}{(1-t^2)^{i+1}},
$$
where $(f_0,\ldots,f_{n-1})$ is the $f$-vector of $K$. This series can be
also expressed in terms of the $h$-vector $(h_0,\ldots,h_n)$
(see~(\ref{hvector})) as
\begin{equation}
\label{poin}
  F\bigl(k(K),t\bigr)=\frac{h_0+h_1t^2+\ldots+h_nt^{2n}}{(1-t^2)^n}.
\end{equation}

On the other hand, the Poincar\'e series of the $k[v_1,\ldots,v_m]$-module
$k(P)$ (or $k(K)$) can be calculated from any free resolution of
$k(P)$. More precisely, the following general theorem holds
(see e.g.,~\cite{St}).

\begin{theorem}
  Let $M$ be a finitely generated graded
  $k[v_1,\ldots,v_m]$-module, $\deg v_i=2$, and there is given a finite free
  resolution of $M$:
  $$
    0\longrightarrow R^{-h}\stackrel{d^{-h}}{\longrightarrow}
    R^{-h+1}\stackrel{d^{-h+1}}{\longrightarrow}\cdots\longrightarrow
    R^{-1}\stackrel{d^{-1}}{\longrightarrow}
    R^{0}\stackrel{d^{0}}{\longrightarrow} M\longrightarrow 0.
  $$
  Suppose that the free $k[v_1,\ldots,v_m]$-modules $R^{-i}$ have their
  generators in dimensions $d_{1i},\ldots,d_{q_ii}$, where
  $q_i=\dim_{k[v_1,\ldots,v_m]}R^{-i}$. Then the Poincar\'e series of $M$
  can be calculated by the following formula:
  $$
    F(M,t)=\frac{\sum_{i=0}^{-h}(-1)^i(t^{d_{1i}}+\ldots+t^{d_{q_ii}})}
    {(1-t^2)^m}.\quad\square
  $$
\end{theorem}
Now let us apply this theorem to the minimal resolution~(\ref{resol}) of
$k(P)=k[v_1,\ldots,v_m]/\!I$. Since all differentials of the
complex~(\ref{timesk}) are trivial, we obtain
\begin{equation}
\label{GS}
  F\bigl(k(P),t\bigr)=(1-t^2)^{-m}\sum_{i=0}^{m-n}(-1)^iF(\T^i,t).
\end{equation}
Combining this with (\ref{algdual}), we get
\begin{multline*}
  F\bigl(k(P),t\bigr)=(1-t^2)^{-m}\sum_{i=0}^{m-n}(-1)^i t^{2m}
  F(\T^{m-n-i},{\textstyle\frac1t})=\\
  =\bigl(1-({\textstyle\frac1t})^2\bigr)^{-m}\cdot
  (-1)^m\sum_{j=0}^{m-n}(-1)^{m-n-j}
  F(\T^j,{\textstyle\frac1t})
  =(-1)^nF\bigl(k(P),{\textstyle\frac1t}\bigr).
\end{multline*}
Substituting here the expressions from the right-hand side of~(\ref{poin})
for $F\bigl(k(P),t\bigr)$ and $F\bigl(k(P),{\frac1t}\bigr)$,
we finally deduce
\begin{equation}
\label{DS}
  h_i=h_{n-i}.
\end{equation}
These are the well-known {\it Dehn--Sommerville equations}~\cite{Br} for
simple (or simplicial) polytopes.

Thus, we see that the algebraic duality~(\ref{algdual}) and the combinatorial
Dehn--Sommerville equations~(\ref{DS}) follow from the Poincar\'e duality
for the manifold $\ZP$. Furthermore, combining~(\ref{poin}) and~(\ref{GS}) we
obtain
\begin{equation}
\label{htor}
  \sum_{i=0}^{m-n}(-1)^iF(\T^i,t)=(1-t^2)^{m-n}h(t^2),
\end{equation}
where $h(t)=\sum_{i=0}^nh_it^i$.

\medskip

We define the subcomplex $\mathcal A$ of the cochain complex
$\bigl[k(P)\otimes\Lambda[u_1,\ldots,u_m],d\bigr]$ from Theorem~\ref{mult}
as follows.  The $k$-module $\mathcal A$ is generated by monomials
$v_{i_1}\ldots v_{i_p}\otimes u_{j_1}\ldots u_{j_q}$ and $1\otimes
u_{j_1}\ldots u_{j_k}$ such that $\{v_{i_1},\ldots,v_{i_p}\}$ spans a
simplex in $K_P$ and $\{i_1,\ldots,i_p\}\cap\{j_1,\ldots,j_q\}=\emptyset$.
It can be easily checked that $d({\mathcal A})\subset{\mathcal A}$ and,
therefore, $\mathcal A$ is a cochain subcomplex. Moreover, $\mathcal A$
inherits the bigraded module structure from
$k(P)\otimes\Lambda[u_1,\ldots,u_m]$ with differential $d$ adding $(1,0)$
to bidegree.

\begin{lemma}
\label{iscoh}
  The cochain complexes
  $\bigl[k(P)\otimes\Lambda[u_1,\ldots,u_m],d\bigr]$ and
  $[{\mathcal A},d]$ have same cohomologies. Hence, the following isomorphism
  of $k$-modules holds:
  $$
    H[{\mathcal A},d]\cong\Tor_{k[v_1,\ldots,v_m]}\bigl(k(P),k\bigr).
  $$
\end{lemma}
\begin{proof}
It is sufficient to prove that any cocycle
$\omega=v^{\alpha_1}_{i_1}\ldots v^{\alpha_p}_{i_p}\otimes
u_{j_1}\ldots u_{j_q}$
from $k(P)\otimes\Lambda[u_1,\ldots,u_m]$ that does
not lie in $\mathcal A$ is a coboundary. To do this we note that if there is
$i_k\in\{i_1,\ldots,i_p\}\cap\{j_1,\ldots,j_q\}$, then
$d\omega$ contains the summand $v^{\alpha_1}_{i_1}\ldots
v^{\alpha_k+1}_{i_k}\ldots v^{\alpha_p}_{i_p}\otimes
u_{j_1}\ldots\widehat{u}_{i_k}\ldots u_{j_q}$, hence, $d\omega\ne0$ --- a
contradiction. Therefore,
$\{i_1,\ldots,i_p\}\cap\{j_1,\ldots,j_q\}=\emptyset$. If $\omega$
contains at least one $v_k$ with degree $\alpha_k>1$, then since
$d\omega=0$, we have $\omega=\pm
d\bigl(v^{\alpha_1}_{i_1}\ldots v^{\alpha_k-1}_{i_k}\ldots
v^{\alpha_p}_{i_p}\otimes u_{i_k}u_{j_1}\ldots u_{j_q}\bigr)$. Now, our
assertion follows from the fact that all non-zero elements of
$k(P)$ of the type $v_{i_1}\ldots v_{i_p}$ correspond to
simplices of $K_P$.
\end{proof}

Now let us introduce the submodules ${\mathcal A}^{*,2r}\subset{\mathcal A}$,
$r=0,\ldots,m$, generated by monomials $v_{i_1}\ldots v_{i_p}\otimes
u_{j_1}\ldots u_{j_q}\in{\mathcal A}$ such that $p+q=r$. Hence, ${\mathcal
A}^{*,2r}$ is the submodule in $\mathcal A$ consisting of all elements of
internal degree $2r$ (i.e. for any $\omega\in{\mathcal A}^{*,2r}$ one has
$\bideg \omega=(*,2r)$; remember that the internal degree corresponds to the
second grading). It is clear that $\sum_{r=0}^{2m}{\mathcal
A}^{*,2r}={\mathcal A}$. Since the differential $d$ does not change the
internal degree, all ${\mathcal A}^{*,2r}$ are subcomplexes of $\mathcal A$.
The cohomology modules of these complexes are exactly $\T^{i,2r}$ and their
dimensions are the bigraded Betti numbers $b^{-i,2r}(\ZP)$. Let us consider
the Euler characteristics of these subcomplexes:
$$
  \chi_r:=\chi({\mathcal A}^{*,2r})=\sum_{q=0}^m(-1)^q\dim_k{\mathcal A}^{-q,2r}
  =\sum_{q=0}^m(-1)^qb^{-q,2r}(\ZP),
$$
and define
\begin{equation}
\label{chii}
  \chi(t)=\sum_{r=0}^m\chi_rt^{2r}.
\end{equation}
Then it follows from Lemma \ref{iscoh} that
\begin{multline*}
  \chi(t)=\sum_{r=0}^m\sum_{q=0}^m(-1)^q\dim_k{\mathcal A}^{-q,2r}t^{2r}=
  \sum_{q=0}^m(-1)^q\sum_{r=0}^m\dim_kH^{-q}[{\mathcal A}^{*,2r}]t^{2r}\\
  =\sum_{q=0}^m(-1)^q\sum_{r=0}^m\dim_k\T^{q,2r}t^{2r}
  =\sum_{q=0}^m(-1)^qF(\T^q,t),
\end{multline*}
where $\T^{q,2r}=H^{-q,2r}\bigl[k(P)\otimes\Lambda[u_1,\ldots,u_m],d\bigr]=
\Tor^{-q,2r}_{k[v_1,\ldots,v_m]}\bigl(k(P),k\bigr)$.
Combining this with formula (\ref{htor}), we get
\begin{equation}
\label{hchi}
  \chi(t)=(1-t^2)^{m-n}h(t^2).
\end{equation}
This formula can be also obtained directly from the definition of $\chi_r$.
Indeed, it can bee easily seen that
\begin{equation}
\label{aqr}
  \dim_k{\mathcal A}^{-q,2r}=f_{r-q-1}\binom{m-r+q}q,\quad
  \chi_r=\sum_{j=0}^m(-1)^{r-j}f_{j-1}\binom{m-j}{r-j},
\end{equation}
(here we set $\binom jk=0$ if $k<0$). Then
\begin{multline}
\label{chidir}
  \chi(t)=\sum_{r=0}^m\chi_rt^{2r}=
  \sum_{r=0}^m\sum_{j=0}^mt^{2j}t^{2(r-j)}(-1)^{r-j}f_{j-1}
  \binom{m-j}{r-j}\\
  =\sum_{j=0}^mf_{j-1}t^{2j}(1-t^2)^{m-j}=
  (1-t^2)^m\sum_{j=0}^nf_{j-1}(t^{-2}-1)^{-j}.
\end{multline}
Further, it follows from~(\ref{hvector}) that
$$
  t^nh(t^{-1})=(t-1)^n\sum_{i=0}^nf_{i-1}(t-1)^{-i}.
$$
Substituting here $t^{-2}$ for $t$ and taking into account~(\ref{chidir}), we
finally obtain
$$
  \frac{\chi(t)}{(1-t^2)^m}=\frac{t^{-2n}h(t^2)}{(t^{-2}-1)^n}=
  \frac{h(t^2)}{(1-t^2)^n},
$$
which is equivalent to (\ref{hchi}).

Formula (\ref{hchi}) allows to express the $h$-vector of a simple polytope
$P^n$ in terms of the bigraded Betti numbers $b^{-q,2r}(\ZP)$ of
the corresponding manifold $\ZP$.

\begin{lemma}
\label{bps}
  The Poincar\'e series
  $F({\mathcal A}^{*,*},\tau,t)=\sum_{r,q}\dim_k{\mathcal A}^{-q,2r}\tau^{-q}t^{2r}$
  of the bigraded module ${\mathcal A}^{*,*}$ is as follows
  $$
    F({\mathcal A}^{*,*},\tau,t)=
    \sum_jf_{j-1}\left(1+\frac{t^2}{\tau}\right)^{m-j}t^{2j}.
  $$
\end{lemma}
\begin{proof}
Using formula (\ref{aqr}), we calculate
\begin{multline*}
  \sum_{r,q}\dim_k{\mathcal A}^{-q,2r}\tau^{-q}t^{2r}=
  \sum_{r,q}f_{r-q-1}\binom{m-r+q}q\tau^{-q}t^{2r}\\
  =\sum_{r,j}f_{j-1}\binom{m-j}{r-j}\tau^{-(r-j)}t^{2r}
  =\sum_jf_{j-1}\left(1+\frac{t^2}{\tau}\right)^{m-j}t^{2j}.
\end{multline*}
\end{proof}

The bigraded Betti numbers $b^{-i,2j}(\ZP)$ can be calculated either by
means of Theorem~\ref{mult} and the results of subsection~4.3 (as we did
before) or by means of the following theorem, which reduces
their calculation to calculating the cohomology of certain subcomplexes
of the simplicial complex $K^{n-1}$ dual to $\partial P^n$.

\begin{theorem}[\rm Hochster, see~\cite{Ho}, \cite{St}]
\label{hoch}
  Let $K$ be a simplicial complex on the vertex set
  $V=\{v_1,\ldots,v_m\}$, and let $k(K)$ be its face ring.
  Then the Poincar\'e series of
  $\T^i=\Tor_{k[v_1,\ldots,v_m]}^{-i}\bigl(k(K),k\bigr)$ is calculated as
  follows
  $$
    F(\T^i,t)=
    \sum_{W\subseteq V}\bigl(\dim_k\tilde{H}_{|W|-i-1}(K_W)\bigr)t^{2|W|},
  $$
  where $K_W$ is the subcomplex of $K$ consisting of all simplices with
  vertices in $W$.\ep
\end{theorem}

However, easy examples show that the calculation based on the above
theorem becomes very involved even for small complexes $K$. It can be
shown also that applying the discussed above result of~\cite{GM}
(see subsection~2.2) to $U(P^n)$ gives the same description of
$H^{*}\bigl(U(P^n)\bigr)$ as that of $H^{*}(\ZP)$ given by the Hochster
theorem. This, of course, conforms with our results from subsection~2.2.

\begin{lemma}
\label{restr}
  For any simple polytope $P$ holds
  $$
    \Tor^{-q,2r}_{k[v_1,\ldots,v_m]}\bigl(k(P),k\bigr)=0
    \qquad\text{for }\;0<r\le q.
  $$
\end{lemma}
\begin{proof}
This can be seen either directly from the construction of
the minimal resolution~(\ref{resol}), or from Theorem~\ref{hoch}.
\end{proof}

\begin{theorem}
\label{lowhom}
  We have
  \begin{enumerate}
  \item $H^1(\ZP)=H^2(\ZP)=0$.
  \item The rank of the third cohomology group of $\ZP$ (i.e. the third
  Betti number $b^3(\ZP)$) equals the number of pairs of vertices of
  the simplicial complex $K^{n-1}$ that are not connected by an edge.
  Hence, if $f_0=m$ is the number of vertices of $K$ and $f_1$ is the
  number of edges, then
  $$
    b^3(\ZP)=\frac {m(m-1)}2-f_1.
  $$
  \end{enumerate}
\end{theorem}
\begin{proof}
It follows from Theorem~\ref{cohomZ} and Lemma~\ref{restr} that
$$
  H^3({\mathcal Z})=\Tor^{-1,4}_{k[v_1,\ldots,v_m]}\bigl(k(P),k\bigr)
  =\T^{1,4}.
$$
By Theorem~\ref{hoch},
$$
  b^{-1,4}(\ZP)=
  \dim_k\T^{1,4}=\sum_{W\subseteq V,|W|=2}\dim_k\tilde{H}_0(K_W).
$$
Now the theorem follows from the fact that $\dim_k\tilde{H}_0(K_W)=0$ if
$K_W$ is a 1-simplex, and $\dim_k\tilde{H}_0(K_W)=1$ if
$K_W$ is a pair of disjoint vertices.
\end{proof}

\begin{remark}
Combining Theorems \ref{cohomZ}, \ref{hoch} and Lemma
\ref{restr} we can also obtain that
$$
  b^4({\mathcal Z})=\dim_k\T^{2,6}=
  \sum_{W\subseteq V,|W|=3}\dim_k\tilde{H}_0(K_W).
$$
\end{remark}

Manifolds $\ZP$ allow to give a nice interpretation not only to the
Dehn--Sommerville equations~(\ref{DS}) but also to a number of other
combinatorial properties of simple polytopes. In particular, using
formula~(\ref{hchi}) one can express the well-known MacMullen
inequalities, the Upper and the Lower Bound Conjectures (see.~\cite{Br}) in
terms of the cohomology of $\ZP$. We review here only two examples.

The first non-trivial MacMullen inequality for a simple polytope $P^n$ can be
written as $h_1\le h_2$ for $n\ge3$. In terms of the $f$-vector this means
that $f_1\ge mn-\binom{n+1}2$. Theorem~\ref{lowhom} shows that
$b^3(\ZP)=\binom m2-f_1$. Hence, we have the following upper bound
for $b^3(\ZP)$:
\begin{equation}
\label{Zbound}
  b^3(\ZP)\le\binom{m-n}2\quad\text{ if }n\ge3.
\end{equation}

The Upper Bound Conjecture for the number of faces of a simple
polytope can be formulated in terms of the $h$-vector as
\begin{equation}
\label{ubc}
  h_i\le\binom{m-n+i-1}i.
\end{equation}
Using the decomposition
$$
  \left(\frac1{1-t^2}\right)^{m-n}=
  \sum_{i=0}^{\infty}\binom{m-n+i-1}it^{2i},
$$
we deduce from~(\ref{hchi}) and~(\ref{ubc}) that
\begin{equation}
\label{chibound}
  \chi(t)\le1, \quad 0\le t<1.
\end{equation}

It would be interesting to obtain a purely topological proof of
inequalities~(\ref{Zbound}) and~(\ref{chibound}).


\newpage

\begin{thebibliography}{GM}

\bibitem[Ad]{Ad}
J.\,F. Adams,
{\it On the non-existence of elements of Hopf invariant one},
Annals of Math. {\bf 72} (1960), no.~1, 20--104.

\bibitem[Ba]{Ba}
V.\,V. Batyrev,
{\it Quantum Cohomology Rings of Toric Manifolds},
Journ\'ees de G\'eometrie Alg\'ebrique d'Orsay (Juillet 1992),
Ast\'erisque {\bf 218}, Soci\'ete Math\'ematique de France, Paris, 1993,
pp.~9--34; available at http://xxx.lanl.gov/find/math.AG/9310004

\bibitem[Br]{Br}
A. Br\o nsted,
{\it An introduction to convex polytopes},
Springer-Verlag, New-York, 1983.

\bibitem[BP]{BP}
V.\,M.~Bukhshtaber and T.\,E.~Panov,
{\it Algebraic topology of manifolds defined by simple polytopes}
(Russian), Uspekhi Mat. Nauk {\bf 53} (1998), no.~3, 195--196;
English transl. in:
Russian Math. Surveys {\bf 53} (1998), no.~3, 623--625.

\bibitem[CE]{CE}
H. Cartan and S. Eilenberg,
{\it Homological algebra},
Princeton Univ. Press, Princeton, N.J., 1956.

\bibitem[Da]{Da}
V.~Danilov,
{\it The geometry of toric varieties},
(Russian), Uspekhi Mat. Nauk {\bf 33} (1978), no.~2, 85--134;
English transl. in:
Russian Math. Surveys {\bf 33} (1978), 97--154.

\bibitem[DJ]{DJ}
M.~Davis and T.~Januszkiewicz,
{\it Convex polytopes, Coxeter orbifolds and torus actions},
Duke Math. Journal {\bf 62} (1991), no.~2, 417--451.

\bibitem[EM]{EM}
S. Eilenberg and J.\,C. Moore,
{\it Homology and fibrations. I},
Comment. Math. Helv. {\bf 40} (1966), 199--236.

\bibitem[Fu]{Fu}
W. Fulton,
{\it Introduction to Toric Varieties},
Princeton Univ. Press, Princeton, N.J., 1993.

\bibitem[GM]{GM}
M. Goresky and R. MacPherson,
{\it Stratified Morse Theory},
Springer-Verlag, Berlin-New York, 1988.

\bibitem[Ho]{Ho}
M. Hochster,
{\it Cohen--Macaulay rings, combinatorics, and simplicial complexes},
in: Ring Theory II (Proc. Second Oklahoma Conference),
B.\,R.~McDonald and R.~Morris, editors, Dekker, New York, 1977,
pp.~171--223.

\bibitem[La]{La}
P.\,S. Landweber,
{\it Homological properties of comodules over $MU_{*}(MU)$ and $BP_{*}(BP)$},
American Journal of Mathematics
{\bf 98} (1976), 591--610.

\bibitem[Ma]{Ma}
S. Maclane,
{\it Homology},
Springer-Verlag, Berlin, 1963.

\bibitem[Se]{Se}
J.--P. Serre,
{\it Alg\`ebre locale-multipliciti\'es},
Lecture Notes in Mathematics {\bf 11}, Springer-Verlag, Berlin, 1965.

\bibitem[Sm]{Sm}
L. Smith,
{\it Homological Algebra and the Eilenberg--Moore Spectral Sequence},
Transactions of American Math. Soc. {\bf 129} (1967), 58--93.

\bibitem[St]{St}
R. Stanley,
{\it Combinatorics and Commutative Algebra},
Progress in Mathematics {\bf 41}, Birkhauser, Boston, 1983.

\end{thebibliography}
\end{document}